\documentclass{article}
\usepackage{amsmath,amsfonts,amssymb,amsthm,pb-diagram,picinpar,graphicx,color}
\usepackage{algorithm}
\usepackage{algorithmic}
\usepackage{colortbl}
\usepackage{enumerate}
\usepackage{hyperref}
\usepackage{amscd}
\usepackage{bm}
\usepackage{pgf}

\usepackage[T1]{fontenc}
\usepackage{tikz}
\usetikzlibrary{cd}
\usepackage{enumitem}
\usepackage{hyperref}
\usepackage{nccmath}
\usepackage{here}
\usepackage{longtable}

\usetikzlibrary{shapes,matrix,calc,arrows,decorations.markings,knots,patterns,intersections,cd}

\def \qed{\hfill $\Box$}

\newcommand{\CC}{\mathbb{C}}
\newcommand{\QQ}{\mathbb{Q}}
\newcommand{\ZZ}{\mathbb{Z}}
\newcommand{\PP}{\mathbb{P}}
\newcommand{\FF}{\mathbb{F}}

\newcommand{\mcC}{\mathcal{C}}

\newcommand{\mcO}{\mathcal{O}}

\newcommand{\mcQ}{\mathcal{Q}}

\newcommand{\fd}{\mathfrak{d}}

\newcommand{\vc}{\boldsymbol {c}}

\newcommand{\I}{\mathop{\mathrm{I}}\nolimits}
\newcommand{\III}{\mathop{\mathrm{III}}\nolimits}
\newcommand{\IV}{\mathop{\mathrm{IV}}\nolimits}

\newcommand{\Div}{\mathop{\mathrm{Div}}\nolimits}

\newcommand{\MW}{\mathop{\mathrm{MW}}\nolimits}
\newcommand{\NS}{\mathop{\mathrm{NS}}\nolimits}

\newcommand{\Supp}{\mathop{\mathrm{Supp}}\nolimits}

\newcommand{\sing}{\mathop{\mathrm{Sing}}\nolimits}
\newcommand{\Red}{\mathop{\mathrm{Red}}\nolimits}

\newcommand{\ord}{\mathop{\mathrm{ord}}\nolimits}

\newcommand{\Kbar}{\overline{K}}

\newtheorem{thm}{Theorem}[section]
\newtheorem{lem}[thm]{Lemma}     
\newtheorem{cor}[thm]{Corollary}
\newtheorem{prop}[thm]{Proposition}

\theoremstyle{definition}
\newtheorem{defin}[thm]{Definition}

\theoremstyle{remark}
\newtheorem{rem}[thm]{Remark}

\title{The geometry of three sections on certain rational elliptic surfaces and Mumford representations}
\date{\today}
\author{Ryosuke MASUYA}

\begin{document}

\maketitle
\begin{abstract}
In this article, we study the geometry of plane curves obtained by three sections and another section given as their sum on certain rational elliptic surfaces.
We make use of Mumford representations of semi-reduced divisors in order to study the geometry of sections.
As a result, we are able to give new proofs for some classical results on singular plane quartics and their bitangent lines.
\\


\end{abstract}

\section*{Introduction}\label{sec:intro}

Let $\varphi:S\rightarrow C$ be an elliptic surface over a smooth projective curve $C$ satisfying the conditions as follows:
{\rm (i)} $\varphi$ is relatively minimal, 
{\rm (ii)} $\varphi$ has a section $O:C\rightarrow S$ and, 
{\rm (iii)} $\varphi$ has at least one singular fiber.
Under these conditions,  the N\'{e}ron-Severi group $\NS(S)$ of $S$ is finitely generated and torsion-free by \cite[Theorem 1.2]{Shio}.
The base field of this article is always the field of complex numbers $\CC$.

Let $E_S$ be the generic fiber of $\varphi$. 
$E_S$ can be regarded as a curve of genus 1 defined over the field of rational functions, $\CC(C)$, of $C$. 
Let $\MW(S)$ be the set of sections of $\varphi$ and let $E_S(\CC(C))$ be the set of $\CC(C)$-rational points of $E_S$. 
It is well-known that we can identify $\MW(S)$ with $E_S(\CC(C))$. 
For a section $s:C\rightarrow S$ we identify $s$ with its image on $S$. 
For $P\in E_S(\CC(C))$, we denote the corresponding section by $s_P$. 
In this article, we also write $O$ for $s_O$. 
For a section $s$, we denote the corresponding rational point of $E_S(\CC(C))$ by $P_s$.
We will denote the sum of two points $P,Q\in E_S$ as divisors by $P+Q$ and the sum as the addition in terms of the group law of $E_S$ by $P\dot{+}Q$.

In this article, we study the geometry of plane quartic curve $\mcQ$ and its bitangent lines by properties of an elliptic surface $S$ and $\MW(S)$, which is obtained by a double cover of $\PP^2$ branched along $\mcQ$.
In order to explain our results and motivation, 
we explain our setting following \cite{Toku10}.

Let $\mcQ$ be a reduced plane quartic curve which is not the union of 4 concurrent lines, let $z_o$ be a smooth point of $\mcQ$ and let $\Lambda_{z_o}$ be the pencil of lines through $z_o$.
Let $q_1:(\PP_{z_o}^2)' \rightarrow \PP^2$ be the blowing-up at $z_o$ and we denote the proper transform of the tangent line $l_{z_o}$ at $z_o$  and the exceptional curve of $q_1$ by $\overline{l}_{z_o}$ and $\Delta_1$.
We next denote the blowing-up at $\overline{l}_{z_o}\cap \Delta_1$ by $q_2:\PP^2_{z_o}\rightarrow (\PP_{z_o}^2)'$.
Put $q_{z_o}=q_1\circ q_2$. 

Let $f':S'_{\mcQ,z_o} \rightarrow \PP^2_{z_o}$ be the double cover with branch locus $\overline{\mcQ}$ and $\Delta$, where $\overline{\mcQ}$ and $\Delta$ is the proper transform of $\mcQ$ and $\Delta_1$ with respect to $q_{z_o}$, respectively.
Let $q$ be a composition of a finite number of blowing-ups so that the $\overline{\mcQ}$ becomes smooth and 
let $\mu: S_{\mcQ,z_o}\rightarrow S'_{\mcQ,z_o} $ be the canonical resolution of $S'_{\mcQ,z_o}$ (see \cite{Hori} for canonical resolutions).
Then we see that $S_{\mcQ,z_o}$ is an elliptic surface satisfying the following properties:

\begin{itemize}
    \item[(i)] The pencil $\Lambda_{z_o}$ of lines through $z_o$ induces a relatively minimal elliptic fibration $\varphi_{\mcQ,z_o}:S_{\mcQ,z_o}\rightarrow \PP^1$.
    \item[(ii)] The preimage of $\Delta$ by $f'\circ \mu$ gives rise to a section $O$ of $\varphi_{\mcQ,z_o}$.
    \item[(iii)] The map $\varphi_{\mcQ,z_o}$ has a singular fiber $F_{z_o}$ containing irreducible components mapped by $f'\circ\mu$ to the exceptional curve of $q_2$ and the proper transform of $\overline{l}_{z_o}$ by $q_2$.
\end{itemize}
Put $\overline{f}_{\mcQ,z_o}=q_{z_o}\circ f'\circ \mu$.

In order to explain our main result, we define {\it line-sections} and {\it weak-bitangent lines} as follows:
\begin{defin}
    For $s\in\MW(S_{\mcQ,z_o})$ and a line $L$, $s$ is said to be a {\it line-section} of $L$ if $\overline{f}_{\mcQ,z_o}(s)$ is a line $L$ in $\PP^2$.
\end{defin}

\begin{defin}
    Let $L$ be a line meeting $\mcQ$ at two points $x_0$ and $y_0$ only.
    If the two local intersection multiplicities $I_{x_0}(\mcQ,L)$ and $I_{y_0}(\mcQ,L)$ are even, $L$ is said to be a {\it weak-bitangent line}.
    If $x_0$ and $y_0$ are smooth points on $\mcQ$ and $x_0\neq y_0$ (resp. $x_0=y_0$), $L$ is said to be a {\it bitangent line } (resp. a {\it $4$-fold tangent line}).
\end{defin}
\begin{rem}
    In Section\,\ref{sec:rep-w-bi}, we see that weak-bitangent lines give rise to line-sections in $\MW(S_{\mcQ,z_o})$.
\end{rem}

The geometry of plane quartic curves and its weak-bitangent lines have been studied by various mathematicians.
For example, Shioda studied a smooth plane quartic curve $\mcQ_1$ and its 28 bitangent lines from the viewpoint of rational points of $E_{S_{\mcQ_1,z_o}}(\CC(t))$ which correspond to line-sections given by the bitangent lines \cite{Shio93}.
Also, Bannai and Tokunaga studied the embedded topology of plane curve arrangement of a certain singular plane quartic curve $\mcQ_2$, weak-bitangent lines and smooth conics by properties of $E_{S_{\mcQ_2,z_o}}(\CC(t))$ and $S_{\mcQ_2,z_o}$ \cite{Ban-Toku15, Ban-Toku17, Ban-Toku21}.
In this article, we are interested in the geometry of the line-sections corresponding to weak-bitangent lines. 
In particular, in this article we assume that $\mcQ$ is a reduced plane quartic curve that satisfies $({\dagger})$ below
\begin{itemize}
	\item[$({\dagger})$]
	an irreducible singular quartic curve, or $\mcQ=\mcC_1+\mcC_2$, where $\mcC_1$ and $\mcC_2$ are two smooth conics meeting transversally. 
\end{itemize}
In this article, by properties of $E_{S_{\mcQ,z_o}}(\CC(t))$ and $S_{\mcQ,z_o}$, we obtain our main result Theorem \ref{thm:main} as follows, which concerns geometry of $\mcQ$ and three weak-bitangent lines.

Our key tool to treat curves is `the Mumford representation of semi-reduced divisor' on an elliptic curve $E_{S_{\mcQ,z_o}}$ over $\CC(t)$, by which we can compute the sum of points on $E_{S_{\mcQ,z_o}}$ defined over $\CC(t)$.
Namely, we can compute the sum of points on $E_{S_{\mcQ,z_o}}$ without considering suitable algebraic extensions of $\CC(t)$ \cite{Ban-Kaw-Mas-Tok}.
Another advantage to utilize Mumford representations is that we are able to calculate the sum of three rational points at one time.
As an application, we study the relation between three line-sections $s_1,s_2$ and $s_3$ and the line-section corresponding to $P_{s_1}\dot{+}P_{s_2}\dot{+}P_{s_3}$.
As a result, we give new proofs for classical results for $\mcQ$ and its weak-bitangent lines.

\begin{thm}\label{thm:main}
    Let $\mcQ$ be a plane curve satisfying $(\dagger)$ and let $z_o$ be a general smooth point of $\mcQ$.
    Let $s_1,s_2$ and $s_3$ be three line-sections of $\MW(S_{\mcQ,z_o})$ and put $P_4=P_{s_1}\dot{+}P_{s_2}\dot{+}P_{s_3}$.
    If the corresponding section $s_{P_4}$ is a line-section and $\overline{f}_{\mcQ,z_o}(s_i)\neq \overline{f}_{\mcQ,z_o}(s_j)$ for all distinct $i,j\in \{1,2,3\}$, then all the points of intersections of $\mcQ$ with the four lines $\overline{f}_{\mcQ,z_o}(s_1), \overline{f}_{\mcQ,z_o}(s_2)$, $\overline{f}_{\mcQ,z_o}(s_3)$ and $\overline{f}_{\mcQ,z_o}(s_{P_4})$ lie on a conic.
\end{thm}

As applications, we give new proofs to the corollaries below of bitangent lines that Harris proved in \cite{Harris}.

\begin{cor}{\rm (\cite[Theorem 3.3]{Harris})}\label{cor:ha3.3}
    Let $C_1$ and $C_2$ on $\PP^2$ be smooth conics meeting transverselly.
    Then the eight points of contact of $C_1+C_2$ with its four bitangent lines all lie on a smooth conic.
\end{cor}
\begin{cor}{\rm (\cite[Theorem 3.4]{Harris})}\label{cor:ha3.4}
    If $\mcQ\subset\PP^2$ is an irreducible quartic with $3$ nodes, then the eight points of contact of $\mcQ$ with its four bitangent lines all lie on a smooth conic.
\end{cor}
\begin{cor}{\rm (\cite[Theorem 3.5]{Harris})}\label{cor:ha3.5}
    An irreducible quartic with an ordinary triple point has four bitangent lines, whose eight points of contact all lie on a smooth conic.
\end{cor}

By Theorem \ref{thm:main}, we also give another theorem which is not given in \cite{Harris}.

\begin{thm}\label{thm:A2+A1}
    Let $\mcQ$ be an irreducible quartic curve with exactly two singularities $x$ and $y$ such that $x$ (resp. $y$) is a simple cusp (resp. a node).
    Then there exist four weak-bitangent lines $L_1,L_2,L_3$ and $L_4$ passing through $x$ and not $y$ and there exist three weak-bitangent lines $M_1,M_2$ and $M_3$ passing through $y$ and not $x$.
    Moreover, when $I_{x}(\mcQ,L_i)=2 \ (i=1,2,3,4)$, for each pair $(L_i,L_j) \ (1\leq i<j\leq 4)$, there exists a unique pair $(M_{a_{ij}},M_{b_{ij}}) \ (1\leq a_{ij}<b_{ij} \leq 3)$ and a smooth conic $C_{ij}$ depending on $(L_i,L_j)$ such that
    \begin{itemize}
        \item[$(\ast)$]
        the six points in $\mcQ\cap(L_i+L_j+M_{a_{ij}}+M_{b_{ij}})$ all lie on $C_{ij}$.
    \end{itemize}
\end{thm}

The orientation of this article is as follows: In Section \ref{sec:elli-sur} we summarize known results of elliptic surfaces.
In Section \ref{sec:rep-w-bi} we give two lists on curves to which sections of $\MW(S_{\mcQ,z_o})$ are mapped by $\overline{f}_{\mcQ,z_o}$. 
In Section \ref{sec:mum-rep} we explain the Mumford representations of semi-reduced divisors which is our key tool.
In Section \ref{sec:appli} we prove the above theorems and corollaries.

\section{Elliptic surfaces}\label{sec:elli-sur}
\subsection{Notation and terminology}
We refer to \cite{Koda, Shio, Mira} for details.
We here define some notation and terminology. 
In this article, an elliptic surface always satisfies the three conditions in the Introduction. 
Let $\varphi:S\rightarrow C$ be an elliptic surface over a smooth projective curve $C$. 
For $v\in C$, we denote the corresponding fiber by $F_v=\varphi^{-1}(v)$. 
We define two subsets, $\sing(\varphi)$ and $\Red(\varphi)$, of $C$  concerning singular fibers as follows:
\begin{align*}
    \sing(\varphi)&:=\{v\in C\mid F_v \ {\rm is \ singular} \}, \\
    \Red(\varphi)&:=\{v\in \sing(\varphi)\mid F_v \ {\rm is \ reducible} \}.
\end{align*}
For $v\in \Red(\varphi)$, the irreducible decomposition of $F_v$ is denoted by
\begin{equation}\label{eq:notation1}
    F_v=\Theta_{v,0}+\sum^{m_v-1}_{i=1}a_{v,i}\Theta_{v,i},
\end{equation}
where $\Theta_{v,0}$ is the unique component with $\Theta_{v,0}\cdot O=1$.
We call $\Theta_{v,0}$ the {\it identity component} of $F_v$. In order to describe the type of singular fibers, we use Kodaira's symbols (\cite{Koda}). Irreducible components of singular fibers are labeled as in \cite{Toku}.

By assumption, for $O\in \MW(S)$, $\MW(S)$ is endowed with an abelian group structure by considering fiberwise addition with $O$ as the zero element.

For $v\in\Red(\varphi)$, we define:
\begin{align}\label{eq:notation2}
    &\begin{aligned}
        \vc(v,D)&:=
        \left[
           \begin{array}{c}
                D\cdot\Theta_{v,1}\\
              \vdots \\
              D\cdot\Theta_{v,m_v-1}
            \end{array}
        \right] \in\ZZ^{m_v-1},\\
        A_v&:=[\Theta_{v,i}\cdot\Theta_{v,j}]_{1\leq i,j \leq m_v-1},\\
        \FF_v&:=[\Theta_{v,1} \cdots \Theta_{v,m_v-1}],
        \end{aligned}
\end{align}
where $D$ is a divisor on $S$.

We explained some notation and terminology of general elliptic surfaces.
 For a rational surface which we consider, we next define notation similarly.
 Let $\mcQ$ be a reduced plane quartic curve satisfying $(\dagger)$ and let $z_o$ be a general smooth point of $\mcQ$.
 We recall the construction of a rational elliptic surface $S_{\mcQ,z_o}$. 
 Since the tangent line $l_{z_o}$ meets $\mcQ$ at $z_o$ and at another two distinct points, $F_{z_o}$ is a singular fiber of type $\I_2$.
 We put $F_{z_o}=\Theta_{\infty,0}+\Theta_{\infty,1}$, where $\Theta_{\infty,0}$ is the identity component of $F_{z_o}$.
 Let $x$ be a singularity of $\mcQ$.
 We denote the singular fiber $\overline{f}_{\mcQ,z_o}^{-1}(x)$ by $F_x$.
 For a singularity $x$ of $\mcQ$ and a divisor $D$ on $S_{\mcQ,z_o}$, the irreducible decomposition of $F_x$, $m_x$ and $\Theta_{x,i}$ are written as in (\ref{eq:notation1}). 
 Moreover, we also define $\vc(x,D)$, $A_x$ and $\FF_x$ as in (\ref{eq:notation2}).

\subsection{Some properties of elliptic surfaces}
We recall some necessary facts.
\begin{thm}{\rm (\cite[Theorem 1.3 ]{Shio})}
    \[
    \overline{\psi}:\NS(S)/T_{\varphi}\simeq E_S(\CC(C))
    \]
    where $T_{\varphi}$ is the subgroup of $\NS(S)$ generated by $O$ and all the irreducible components of fibers. 
\end{thm}
    
Given a divisor $D$ on $S$, we denote $\overline{\psi}(D \ \mathrm{mod}\,T_{\varphi} )$ by $P_D$. 

\begin{thm}{\rm \cite[Lemma 5.1]{Shio}}\label{lem:shio-form}
    For $D\in\Div(S)$, there is a unique section $s(D)$ such that
    \begin{equation*}
        D\approx s(D)+(d-1)O+nF+\sum_{v\in\Red(\varphi)}\FF_vA_v^{-1}\vc(v,D-s(D))
    \end{equation*}
    where $\approx$ and $\chi(\mcO_S)$ denote algebraic equivalence and the characteristic of $\mcO_S$, $d=D\cdot F$ and $n=(d-1)\chi(\mcO_S)+O\cdot(D-s(D))$.
\end{thm}
Note that $A_v^{-1}\vc(v,D-s(D))\in\ZZ^{m_v-1}$ by Theorem\,\ref{lem:shio-form}, while entries of $A_v^{-1}$ are not necessary integers.

\begin{rem}
    For $D\in\Div(S)$, we have $s(D)=s_{P_D}$, where $s(D)$ is the section for $D$ in Theorem\,\ref{lem:shio-form}. 
\end{rem}

\begin{rem}\label{rem:become-line}
    Let $\mcQ$ be a reduced plane quartic curve satisfying $(\dagger)$. 
    Let $P_1,\ldots ,P_n$ be generators of $E_{S_{\mcQ,z_o}}(\CC(t))$ and $Q=\sum_i c_iP_i\in E_{S_{\mcQ,z_o}}(\CC(t)) \ (c_i\in\ZZ)$.
    By Theorem \ref{lem:shio-form} $(\vc(\infty,\sum_ic_is_{P_i})-\vc(\infty,s_{Q}))/2\in\ZZ$ holds.
    Hence, we have
    \[
        s_Q\cdot \Theta_{\infty,1}=
        \left\{
            \begin{array}{cc}
                1& (\sum_ic_is_{P_i}) \cdot \Theta_{\infty,1} {\rm :odd}\\
                0& {\rm otherwise}
            \end{array}
        \right. .
    \]
\end{rem}

Let $\phi_o:\Div(S)\rightarrow \NS(S)_{\QQ}$ and $\phi:E_S(\CC(C))\rightarrow \NS(S)_{\QQ}$ be the homomorphisms given in \cite{Ban-Toku21} and \cite{Shio} respectively.
In \cite{Shio}, Shioda defined a structure on $E_S(\CC(C))$ called the {\it height pairing} denoted by $\langle -,-\rangle$. 
We refer to \cite{Shio} for details. 
Then we have

\begin{thm}[\cite{Shio}]\label{thm:hiegh-pair}
    Let $s_1, s_2\in\MW(S)$. 
    The height pairing $\langle P_{s_1},P_{s_2}\rangle$ is given by
    \[
        \langle P_{s_1},P_{s_2}\rangle = \chi(\mcO_S)-s_1\cdot s_2+s_1\cdot O+s_2\cdot O-\sum_{v\in\Red(\phi)}\mathrm{contr}_v(s_1,s_2),
    \]
    where for divisors $D_1$ and $D_2$ on $S_{\mcQ,z_o}$, $\mathrm{contr}_v(D_1,D_2)$ is given by
    \[
        \mathrm{contr}_v(D_1,D_2)= {}^t \vc(v,D_1)(-A_v)^{-1}\vc(v,D_2).
    \]
\end{thm}

\begin{rem}
    Let $\mcQ$ be a reduced plane quartic curve satisfying $(\dagger)$ and let $z_o$ be a general smooth point of $\mcQ$.
    For a singularity $x$ of $\mcQ$ and divisors $D_1$ and $D_2$ on $S_{\mcQ,z_o}$,we put 
    \begin{align*}
        \mathrm{contr}_x(D_1,D_2)&={}^t\vc(x,D_1)(-A_x^{-1})\vc(x,D_2).
    \end{align*}
\end{rem}

\section{Line-sections arising from weak-bitangents and bitangents and their description in $\MW(S_{\mcQ,z_o})$}\label{sec:rep-w-bi}
In Section \ref{sec:rep-w-bi}, we assume that a plane curve $\mcQ$ satisfies $({\dagger})$ and $z_o$ is a general smooth point on $\mcQ$.
For some rational points $P\in E_{S_{\mcQ,z_o}}(\CC(t))$, we consider the images $\overline{f}_{\mcQ,z_o}(s_P)$.
In particular, in Section \ref{sec:line-sec}, we study its images in the case when $s_P$ are line-sections of weak-bitangent lines.
For this purpose, we give the group structures of $E_{S_{\mcQ,z_o}}(\CC(t))$ in Section \ref{subsec:gene-MW}.
In Section \ref{subsec:node-cusp} and \ref{sec:geo-gene}, we consider the images $\overline{f}_{\mcQ,z_o}(s_P)$, where $P$ are generators given in Section \ref{subsec:gene-MW}.

\subsection{The group structures with height pairing of $E_{S_{\mcQ,z_o}}(\CC(t))$ and propeties of sections}\label{subsec:gene-MW}

In \cite{Toku10}, for all irreducible quartic curves $\mcQ$ with singularities and a general smooth point, $z_o$, on $\mcQ$, Tokunaga gave the group structure with height pairing of $E_{S_{\mcQ,z_o}}(\CC(t))$ based \cite{Ogu-Shio} and, in \cite{Ban-Kaw-Mas-Tok}, for $C_1$ and $C_2$ meeting transversally, the authors gave that of $E_{S_{\mcC_1+\mcC_2,z_o}}(\CC(t))$. 

For all quartic curves $\mcQ$ that satisfy $(\dagger)$, we give the list of the group structure with height pairing (Table \ref{table:st-MW}), following \cite{Ogu-Shio, Toku10, Ban-Kaw-Mas-Tok}. 
Before we go on to give the list Table \ref{table:st-MW}, we need to introduce some notation.
\begin{itemize}
	\item 
	The column labeled $\Xi_{\mcQ}$ indicates the types of singularities of $\mcQ$.
	 We use the notation in \cite[pp.81-82]{BHPV} in order to describe the types of singularities. 
	\item 
	We write  the subgroup of $\NS(S)$ generated by the set $\{\Theta_{v,1},\ldots,\Theta_{v,m_v-1} \}_{v\in\Red(\varphi_{\mcQ,z_o})}$ in the column labeled $R_{\varphi_{\mcQ,z_o}}$.
	Note that the group is isomorphic to a direct sum of root lattices of A-D-E type with respect to the intersection multiplicity, and we describe it as the direct sum.
	\item 
	The column labeled $E_{S_{\mcQ,z_o}}(\CC(t))$ indicate the group structure with the height pairing of $E_{S_{\mcQ,z_o}}(\CC(t))$.
\end{itemize}
\begin{longtable}{|c|c|c|c|}
	\caption{}
	\label{table:st-MW}
	\endhead
	\hline
	No. & $\Xi_\mcQ$ & $R_{\varphi_{\mcQ,z_o}}$ & $E_{S_{\mcQ,z_o}}(\CC(t))$ \rule[-1.5mm]{0mm}{2mm}\\

	\hline
	1&$A_6$ &			$A_6\oplus A_1$ & 									$\langle1/14\rangle$\\
	\hline	
	2&$E_6$ & 			$E_6\oplus A_1$ &									$\langle 1/6 \rangle$\\
	\hline
	3&$A_5$ &			$A_5\oplus A_1$ &									$A_1^*\oplus \langle1/6\rangle$\\
	\hline
	4&$D_5$ &			$D_5\oplus A_1$ &									$A_1^*\oplus \langle 1/4\rangle$\\
	\hline
	5&$D_4$ &			$D_4\oplus A_1$ &									$(A_1^*)^{\oplus 3}$\\
	\hline
	6&$A_4+A_2$&		$A_4\oplus A_2\oplus A_1$&					$\langle1/30\rangle$\\
	\hline
	7&$A_4+A_1$&		$A_4\oplus A_1^{\oplus 2}$&	
	$\frac{1}{10}
	\left[
		\begin{array}{cc}
			2&1\\
			1&3
		\end{array}
	\right]
	$\\
\hline
8&$4A_1$&			$A_1^{\oplus 5}$ &									$(A_1^*)^{\oplus3}\oplus\ZZ/2\ZZ$\\
\hline
9&$A_3+A_2$&		$A_3\oplus A_2\oplus A_1$&					$A_1^*\oplus \langle1/12\rangle$\\
\hline
10&$A_3+A_1$&		$A_3\oplus A_1^{\oplus 2}$&					$(A_1^*)^{\oplus 2}\oplus \langle1/4\rangle$\\
\hline
11&$3A_2$ &			$A_2^{\oplus 3}\oplus A_1$&					$\langle1/6\rangle\oplus \ZZ/3\ZZ$\\
\hline
12&$2A_2+A_1$&	$A_2^{\oplus 2}\oplus A_1^{\oplus 2}$&	$\langle1/6\rangle^{\oplus 2}$\\
\hline
13&$A_2+2A_1$&	$A_2\oplus A_1^{\oplus 3}$&					
$A_1^*\oplus\frac{1}{6}
\left[
\begin{array}{cc}
2&1\\
1&2
\end{array}
\right]
$
\\
\hline
14&$3A_1$&			$A_1^{\oplus 4}$&										$(A_1^*)^{\oplus 4}$\\
\hline
15&$A_4$&				$A_4\oplus A_1$&										
$\frac{1}{10}
\left[
\begin{array}{ccc}
3&1&-1\\
1&7&3\\
-1&3&7
\end{array}
\right]
$
\\
\hline
16&$A_3$ &			$A_3\oplus A_1$	 &										$A_3^*\oplus A_1^*$\\
\hline
17&$2A_2$&			$A_2^{\oplus 2}\oplus A_1$&					$A_2^*\oplus \langle 1/6\rangle$\\
\hline
18&$A_2+A_1$&		$A_2\oplus A_1^{\oplus 2}$&					
$\frac{1}{6}
\left[
\begin{array}{cccc}
2&1&0&-1\\
1&5&3&1\\
0&3&6&3\\
-1&1&3&5
\end{array}
\right]
$\\
\hline
19&$2A_1$&			$A_1^{\oplus 3}$ &									$D_4^*\oplus A_1^*$\\
\hline
20&$A_2$ &		$A_2\oplus A_1$ &		$A_5^*$\\
\hline
21&$A_1$ &		$A_1^{\oplus 2}$ &		$D_6^*$\\
\hline
\end{longtable}

Let $I$ be the matrix corresponding with the structure given in Table \ref{table:st-MW}.
We choose generators $P_1,...,P_n, P_{\tau}\in E_{S_{\mcQ,z_o}}(\CC(t))$ such that
$I=[\langle P_i,P_j\rangle]_{1\leq i,j\leq n}$ and $\ P_{\tau}$ is a torsion of $E_{S_{\mcQ,z_o}}(\CC(t))$.
In Section \ref{sec:geo-gene}, we will study the images $\overline{f}_{\mcQ,z_o}(s_{P_i})$ (Table \ref{table:geo-gene}).
We first consider the case when $\mcQ$ has singularities whose types are $A_1$ and $A_2$.
Note that this case is the case considered in Theorem \ref{thm:A2+A1}.
We can also consider other cases in Table \ref{table:geo-gene}.
For this propose, we introduce the following lemma:

\begin{lem}{\rm (\cite[Lemma 9]{Ban-Toku17})}\label{lem:line-sec}
	Let $s\in\MW(S_{\mcQ,z_o})$ be a section such that $s\cdot O=0$ and $s\cdot \Theta_{\infty,1}=1$. Then $\overline{f}_{\mcQ,z_o}(s)$ is a line $L_s$ such that
	\begin{itemize} 
		\item[$\mathrm{(i)}$]
			$I_x(\mcQ,L_s)$ is even for all $x\in\mcQ$, and 
		\item[$\mathrm{(ii)}$]
			$z_o\notin L_s$.
	\end{itemize}
	Conversely, any line satisfying the two conditions $\mathrm{(i)}$ and $\mathrm{(ii)}$ gives rise to line-sections $s_{L^{\pm}}$ such that $s_{L^{\pm}}\cdot O=0$ and $s_{L^{\pm}}\cdot \Theta_{\infty,1}=1$.
\end{lem}
In the above condition, if $s\cdot\Theta_{\infty,1}=0$ and $s\cdot O=0$ then $\overline{f}_{\mcQ,z_o}(s)$ is a smooth conic that is tangent at $z_o$. 
For details, see \cite{Toku10}.
By definition of weak-bitangent lines and Lemma \ref{lem:line-sec}, weak-bitangent lines give rise to line-sections of $S_{\mcQ,z_o}$.

\subsection{Descriptions of the images of the corresponding sections to the generators of $E_{S_{\mcQ,z_o}}(\CC(t))$ by $\overline{f}_{\mcQ,z_o}$ in the case when $\mcQ$ is No.18}\label{subsec:node-cusp}
Let $x_1$ and $x_2$ be singularities of $\mcQ$ such that the types of $x_1$ and $x_2$ are $A_1$ and $A_2$, respectively. 
By Table \ref{table:st-MW}, the lattice structure of $E_{S_{\mcQ,z_o}}(\CC(t))$ is given by
\[
	I=\frac{1}{6}\left[
		\begin{array}{cccc}
			2&1&0&-1\\
			1&5&3&1\\
			0&3&6&3\\
			-1&1&3&5
		\end{array}
	\right].
\]
Let $P_1,P_2,P_3$ and $P_4$ be generators of $E_{S_{\mcQ,z_o}}(\CC(t))$ such that $I=[\langle P_i,P_j \rangle]$.
Note that the type of $F_{x_1}$ is $\I_2$ or $\III$, that of $F_{x_2}$ is $\I_3$ or $\IV$ and that of $F_{z_o}$ is $\I_2$. 
For $P$ and $Q$ of $E_{S_{\mcQ,z_o}}(\CC(t))$, we have the values of $\mathrm{contr}_{\bullet}(s_P,s_Q)$ as follows:
\begin{align}\label{eq:contr}
	&\begin{aligned}
		\mathrm{contr}_{x_1}(s_P,s_Q)&=
		\left\{
			\begin{array}{cl}
				1/2 & s_P\cdot \Theta_{x_1,1}=s_Q\cdot \Theta_{x_1,1}=1\\
				0 & {\rm otherwise}
			\end{array}
		\right.\\
		\mathrm{contr}_{x_2}(s_P,s_Q)&=
		\left\{
			\begin{array}{cl}
				2/3 & s_P\cdot \Theta_{x_2,i}=s_Q\cdot \Theta_{x_2,i}=1 \ {\rm and} \ i \in \{1,2\}  \\
				1/3 & s_P\cdot \Theta_{x_2,i}=s_Q\cdot \Theta_{x_2,j}=1 \ {\rm and} \ \{ i, j \}=\{1,2\}    \\
				0 & {\rm otherwise}
			\end{array}
		\right.\\
		\mathrm{contr}_{z_o}(s_P,s_Q)&=
		\left\{
			\begin{array}{cl}
				1/2 & s_P\cdot \Theta_{\infty,1}=s_Q\cdot \Theta_{\infty,1}=1\\
				0 & {\rm otherwise}
			\end{array}
		\right.
	\end{aligned}
\end{align}
For $i=1,2,3,4$, We estimate the pairings $\langle P_i,P_i \rangle$:
\begin{align*}
	\langle P_i,P_i \rangle&=2+2s_{P_i}\cdot O-\sum_{y\in J}\mathrm{contr}_y(s_{P_i}\cdot s_{P_i})\qquad (J:=\{ x_1,x_2,z_o \})\\
	&\geq 2+2s_{P_i}\cdot O-(1/2+2/3+1/2)\\
	&=1/3+2s_{P_i}\cdot O 
\end{align*}
Observing that the diagonal components of $I$ are $1/3$, $5/6$ and $1$, we have $s_{P_i}\cdot O=0 \ (i=1,2,3,4)$.
Hence, $\overline{f}_{\mcQ,z_o}(s_{P_i})$ is a line or a conic that is tangent to $\mcQ$ at $z_o$.

We will prove that $s_{P_1}$ is a line-section.
Since the value of $\langle P_1,P_1\rangle$ is $1/3$, we obtain 
\[
	\sum_{y\in J}\mathrm{contr}_y(s_{P_1}\cdot s_{P_1})=5/3.
\]
Therefore, from possible values of $\mathrm{contr}_{y}(s_{P_1},s_{P_1})$ in (\ref{eq:contr}), we obtain 
\begin{align*}
	\mathrm{contr}_{x_1}(s_{P_1},s_{P_1})&=\mathrm{contr}_{z_o}(s_{P_1},s_{P_1})=1/2 \ {\rm and} \\
	\mathrm{contr}_{x_2}(s_{P_1},s_{P_1})&=2/3. 
\end{align*}

By (\ref{eq:contr}), we have $s_{P_1}\cdot\Theta_{x_1,1}=s_{P_1}\cdot \Theta_{\infty,1}=1$ and $s_{P_1}\cdot \Theta_{x_1,i}=1$, where $i=1$ or $2$.
By Lemma \ref{lem:line-sec}, $\overline{f}_{\mcQ,z_o}(s_{P_1})$ is a weak-bitanget line passing through $x_1$ and $x_2$.
We may assume $s_{P_1}\cdot\Theta_{x_2,1}=1$.
Similarly, we find that $s_{P_3}$ is a line-section meeting $\Theta_{x_1,1}$, $\Theta_{x_2,0}$ and $\Theta_{\infty,1}$. 
The image $\overline{f}_{\mcQ,z_o}(s_{P_3})$ is a weak-bitangent line passing through $x_1$ and not $x_2$.
We need to find the images $\overline{f}_{\mcQ,z_o}(s_{P_2})$ and $\overline{f}_{\mcQ,z_o}(s_{P_4})$.
By the assumption, 
\[
	\sum_{y\in J}\mathrm{contr}_y(s_{P_2}\cdot s_{P_2})=\sum_{y\in J}\mathrm{contr}_y(s_{P_4}\cdot s_{P_4})=7/6.
\]
We infer 
\[
	\left[
		\begin{array}{c}
			\mathrm{contr}_{x_1}(s_{P_i}\cdot s_{P_i})\\
			\mathrm{contr}_{x_2}(s_{P_i}\cdot s_{P_i})\\
			\mathrm{contr}_{z_o}(s_{P_i}\cdot s_{P_i})
		\end{array}
	\right]
	=
	\left[
		\begin{array}{c}
			1/2\\
			2/3\\
			0
		\end{array}
	\right]
	\ {\rm or} \
	\left[
		\begin{array}{c}
			0\\
			2/3\\
			1/2
		\end{array}
	\right]
 	\ i=2,4
\]
by the possible values of $\mathrm{contr}_{\bullet}(-,-)$ in (\ref{eq:contr}).
In the former case, $\overline{f}_{\mcQ,z_o}(s_{P_i})$ is a smooth conic passing through $x_1$ and $x_2$.
In the other case, $\overline{f}_{\mcQ,z_o}(s_{P_i})$ is a weak-bitangent passing through $x_2$.
In order to find the values of $\mathrm{contr}_{\bullet}(s_{P_i},s_{P_i}) \ (i=2,4)$, we consider $\langle P_1,P_2\rangle$ and $\langle P_1,P_4\rangle$.

By our choice of generators $P_1,P_2,P_3$ and $P_4$, we have
\begin{align*}
	1/6&=\langle P_1,P_2\rangle=1-s_{P_1}\cdot s_{P_2}-\sum_{y\in J}\mathrm{contr}_y(s_{P_1},s_{P_2}) \ {\rm and}\\
	-1/6&=\langle P_1,P_4\rangle=1-s_{P_1}\cdot s_{P_4}-\sum_{y\in J}\mathrm{contr}_y(s_{P_1},s_{P_4}).
\end{align*}
Now, by (\ref{eq:contr}), the possibilities of $\sum_{y\in J}\mathrm{contr}_y(s_{P_1},s_{P_2})$ are
\[
	0, \ 1/3, \ 1/2, \ 2/3, \ 5/6, \ 1, \ 7/6, \ 4/3 \ {\rm or} \ 5/3.
\]
Since $s_{P_1}\cdot s_{P_2}$ and $s_{P_1}\cdot s_{P_4}$ are integers, we have
\[
	\sum_{y\in J}\mathrm{contr}_y(s_{P_1},s_{P_2})=5/6\ {\rm and} \ \sum_{y\in J}\mathrm{contr}_y(s_{P_1},s_{P_4})=7/6.
\]
The first implies
\begin{align}\label{eq:1}
	\left[
		\begin{array}{c}
			\mathrm{contr}_{x_1}(s_{P_1}\cdot s_{P_2})\\
			\mathrm{contr}_{x_2}(s_{P_1}\cdot s_{P_2})\\
			\mathrm{contr}_{z_o}(s_{P_1}\cdot s_{P_2})
		\end{array}
	\right]
	=
	\left[
		\begin{array}{c}
			1/2\\
			1/3\\
			0
		\end{array}
	\right]
 	\ {\rm or} \  
	\left[
		\begin{array}{c}
			0\\
			1/3\\
			1/2
		\end{array}
	\right].
\end{align}
On the other hand, by the definition of $\mathrm{contr}_{\bullet}(-,-)$, $\mathrm{contr}_{y}(s_{P_1},s_{P_2})$ are given by
\begin{align}\label{eq:2}
	\begin{aligned}
		\mathrm{contr}_{x_1}(s_{P_1},s_{P_2})&=\vc(x_1,s_{P_1})(-A_{x_1}^{-1})\vc(x_1,s_{P_2})\\
		&=s_{P_2}\cdot\Theta_{x_1,1}/2,\\
		\mathrm{contr}_{x_2}(s_{P_1},s_{P_2})&={}^t\vc(x_2,s_{P_1})(-A_{x_2}^{-1})\vc(x_2,s_{P_2})\\
		&=[1 \ 0]
		\left[
			\begin{array}{cc}
				2/3 & 1/3\\
				1/3 & 2/3
			\end{array}
		\right]
		\left[
			\begin{array}{c}
				s_{P_2}\cdot\Theta_{x_2,1}\\
				s_{P_2}\cdot\Theta_{x_2,2}
			\end{array}
		\right],\\
		\mathrm{contr}_{z_o}(s_{P_1},s_{P_2})&=s_{P_2}\cdot\Theta_{\infty,1}/2.
	\end{aligned}
\end{align}
From (\ref{eq:1}) and (\ref{eq:2}), we have $(\vc(x_1,s_{P_2}) , \vc(z_o,s_{P_2}))=(1,0) \ {\rm or} \ (0,1)$ and ${}^t\vc(x_2,s_{P_2})=[0 \ 1]$.
Similarly, $(\vc(x_1,s_{P_4}) , \vc(z_o,s_{P_4}))=(1,0) \ {\rm or} \ (0,1)$ and ${}^t\vc(x_2,s_{P_4})=[1 \ 0]$ hold.

We find the following, for $i=2,4$: 
\begin{itemize}
	\setlength{\leftskip}{-0.5cm}
	\item
	If $(\vc(x_1,s_{P_i}),\vc(z_o,s_{P_i}))=(1,0)$, then $\overline{f}_{\mcQ,z_o}(s_{P_i})$ is a smooth conic passing through $x_1$ and $x_2$. 
	\item
	If $(\vc(x_1,s_{P_i}),\vc(z_o,s_{P_i}))=(0,1)$, then $\overline{f}_{\mcQ,z_o}(s_{P_i})$ is a weak-bitangent passing through $x_2$ and not $x_1$.
\end{itemize}
If one of $s_{P_2}$ and $s_{P_4}$ is a line-section and the other is not, 
then it follows that $\mathrm{contr}_{x_1}(s_{P_2},s_{P_4})=\mathrm{contr}_{z_o}(s_{P_2},s_{P_4})=0$.

On the other hand, we have 
\begin{align*}
	\langle P_2,P_4 \rangle&=1-s_{P_2}\cdot s_{P_4}-\sum_{y\in J}\mathrm{contr}_{y}(s_{P_2},s_{P_4})\\
	&=2/3-s_{P_2}\cdot s_{P_4}.
\end{align*}
But $\langle P_2,P_4 \rangle=1/6$, which is a contradiction.

Hence, $s_{P_2}$ and $s_{P_4}$ are both line-sections or not line-sections. 
Here, we assume that $s_{P_2}$ and $s_{P_4}$ are both not line-sections. 
Let $Q_1=P_1, \ Q_2=\dot{-}P_2\dot{+}P_3, \ Q_3=P_3$ and $Q_4=P_3\dot{-}P_4$. 
We can verify that the matrix $[\langle Q_i,Q_j \rangle]$ equals to $[\langle P_i,P_j \rangle]$.
We can take new generators $Q_1,Q_2,Q_3$ and $Q_4$ of $E_{S_{\mcQ,z_o}}(\CC(t))$. 
Then $s_{Q_2}$ and $s_{Q_4}$ are line-sections because $s_{Q_2}\cdot \Theta_{\infty,1}=s_{Q_4}\cdot \Theta_{\infty,1}=1$.
Therefore we can assume $s_{P_2}$ and $s_{P_4}$ are line-sections.

\subsection{Descriptions of the images of the corresponding sections to the generators of $E_{S_{\mcQ,z_o}}(\CC(t))$ by $\overline{f}_{\mcQ,z_o}$ in the case when $\mcQ$ satisfy $(\dagger)$}\label{sec:geo-gene}

For all $\mcQ$ satisfying $(\dagger)$, we obtain Table \ref{table:geo-gene} by computation similar to the above case.
In Table \ref{table:geo-gene}, we choose generators $P_1,...,P_n, P_{\tau}\in E_{S_{\mcQ,z_o}}(\CC(t))$ such that
$I=[\langle P_i,P_j\rangle]_{1\leq i,j\leq n}$ and $\ P_{\tau}$ is a torsion of $E_{S_{\mcQ,z_o}}(\CC(t))$, where the matrix $I$ corresponds to the group structure with the height pairing given in Table \ref{table:st-MW}.
In the above case, we made $s_{P_2}$ and $s_{P_4}$ line-sections.
Picking a special generators, we can obtain Table \ref{table:geo-gene}.
We explain some notations used in the table.
\begin{itemize}
	\item 
	In the column labeled $\sing(\mcQ)$, we write the all pairs $(x,S_x)$, where $x$ is a singularity of $\mcQ$ and $S_x$ is its type of singularity.
	\item
	The column labeled $E_{S_{\mcQ,z_o}}(\CC(t))$ is the same as those in Table \ref{table:st-MW}.
	\item
	Let $P_1,\ldots,P_n,P_{\tau}$ be generators of  $E_{S_{\mcQ,z_o}}(\CC(t))$ as the group structure with the height pairing given by Table \ref{table:st-MW}.
	The column labeled {\rm COG} will indicate the image $\overline{f}_{\mcQ,z_o}(s_{P_i})$ by using the following notation.
	\begin{itemize}
		\item 
		$BL$: a bitangent line to $\mcQ$ or a 4-fold tangent line at one smooth point.
		\item
		$WL(x)$: a weak-bitangent line meeting $\mcQ$ at a singularity $x$ of $\mcQ$ such that  $I_x(\mcQ,WL(x))=4$.  
		For example, when the type of $x$ is $A_1$, $WL(x)$ is tangent to one of two branches at $x$ of $\mcQ$ with multiplicity 3.
		Also, when the type of $x$ is $D_4$, $WL(x)$ is tangent to one of three branches at $x$ of $\mcQ$ with multiplicity 2 and the other two branches with multiplicities 1.  
 		\item
		$WL(x,y)$: a weak-bitangent line meeting $\mcQ$ at two distinct points $x$ and $y$ of $\mcQ$ such that $I_x(\mcQ,WL(x,y))=I_y(\mcQ,WL(x,y))=2$.
		\item 
		$\eta$: a smooth point on $\mcQ$.
		\item
		$CT(x_1,\ldots, x_l;y_1,\ldots,y_m)$ : a conic such that satisfying
		\begin{itemize}
			\item[{\rm (a)}] 
			$x_1,...,x_l$ and $y_1,...,y_m$ are all singularities of $\mcQ$ through which the conic pass,
			\item[{\rm (b)}]
			$I_{x_i}(CT(x_1,\ldots, x_l;y_1,\ldots,y_m),\mcQ)=\mathrm{mult}_{x_i}(\mcQ) \ (i=1,\ldots, l)$ and
			\item[{\rm (c)}] 
			$I_{y_j}(CT(x_1,\ldots, x_l;y_1,\ldots,y_m),\mcQ)>\mathrm{mult}_{y_j}(\mcQ) \ (j=1,\ldots, m)$,
		\end{itemize}
		where we denote the multiplicity of $\mcQ$ at a point $z$ by $\mathrm{mult}_{z}(\mcQ)$.
		When there are no singularities $x_1,\ldots,x_l$ satisfying {\rm (b)}, we write $CT(-;y_1,\ldots,y_m)$.
		We write $CT(x_1,\ldots,x_l;-)$ similarly.
	\end{itemize}
	Let $P\in E_{S_{\mcQ,z_o}}(\CC(t))$. If $\overline{f}_{\mcQ,z_o}(s_P)$ is a curve $C_P$ of the above curves, we write $PC_P$ in the last column.
	Otherwise, we write the values $s_P\cdot O$, $s_P\cdot\Theta_{x,i}$ and $s_P\cdot \Theta_{\infty,1}$, where  $x$ is a singularity of $\mcQ$ and $i=1,\ldots,m_x-1$.
\end{itemize}
\begin{longtable}{|c|c|c|c|}
\caption{}
\label{table:geo-gene}
\endhead
\hline No.&
$\sing(\mcQ)$&$E_{S_{\mcQ,z_o}}(\CC(t))$ & {\rm COG} \rule[-1.5mm]{0mm}{2mm}\\
\hline 1&
$(x,A_6)$&
$\langle1/14\rangle$&
$P_1WL(x)$
\\
\hline	2&
$(x,E_6)$
&		$\langle 1/6 \rangle$ &
$P_1WL(x)$
\\
\hline 3&
$(x,A_5)$	
&		$A_1^*\oplus \langle1/6\rangle$&
\begin{tabular}{c}
$P_1CT(x)$, $P_2WL(x)$
\end{tabular}
\\
\hline 4&
$(x,D_5)$	
&		$A_1^*\oplus \langle 1/4\rangle$&
$P_1WL(x)$, $P_2WL(x)$
\\
\hline 5&
$(x,D_4)$	
&			$(A_1^*)^{\oplus 3}$&
$P_iWL(x)$  $(i=1,2,3)$
\\
\hline 6&
\begin{tabular}{c}
$(x,A_4)$\\
$(y,A_2)$
\end{tabular}
	&	$\langle1/30\rangle$&
$P_1WL(x,y)$
\\
\hline 7&
\begin{tabular}{c}
$(x,A_4)$\\
$(y,A_1)$
\end{tabular}
		&
$\frac{1}{10}
\left[
\begin{array}{cc}
2&1\\
1&3
\end{array}
\right]
$&
\begin{tabular}{c}
$P_1WL(x,y)$,
 $P_2WL(x)$
\end{tabular}

\rule[-4mm]{0mm}{9.5mm}\\
\hline 8&
\begin{tabular}{c}
$(x,A_1)$\\
$(y,A_1)$\\
$(z,A_1)$\\
$(w,A_1)$
\end{tabular}
			&$(A_1^*)^{\oplus3}\oplus\ZZ/2\ZZ$&

\begin{tabular}{l}
$P_1WL(x,y)$, $P_2L(y,z)$,\\ 
$P_3WL(x,z)$, $P_4CT(x,y,z,w;-)$
\end{tabular}
\\
\hline 9&
\begin{tabular}{c}
$(x,A_3)$\\
$(y,A_2)$
\end{tabular}
			&$A_1^*\oplus \langle1/12\rangle$&
$P_1WL(x)$, $P_2WL(x,y)$
			\\
\hline 10&
\begin{tabular}{c}
$(x,A_3)$\\
$(y,A_1)$
\end{tabular}			&$(A_1^*)^{\oplus 2}\oplus \langle1/4\rangle$&
$P_1CT(y;x)$, $P_2WL(x)$, $P_3WL(x,y)$
\\
\hline 11&
\begin{tabular}{c}
$(x,A_2)$\\
$(y,A_2)$\\
$(z,A_2)$
\end{tabular}
		&$\langle1/6\rangle\oplus \ZZ/3\ZZ$&
$P_1WL(x,y)$, $P_2CT(x,y,z;-)$
\\
\hline 12&
\begin{tabular}{c}
$(x,A_2)$\\
$(y,A_2)$\\
$(z,A_1)$
\end{tabular}
		&$\langle1/6\rangle^{\oplus 2}$&
$P_1WL(x,y)$, $P_2CT(x,y,z;-)$
\\

\hline 13&
\begin{tabular}{c}
$(x,A_2)$\\
$(y,A_1)$\\
$(z,A_1)$
\end{tabular}			&	
$A_1^*\oplus\frac{1}{6}
\left[
\begin{array}{cc}
2&1\\
1&2
\end{array}
\right]
$&
$P_1WL(y,z)$, $P_2WL(x,y)$, $P_3WL(x,z)$
\\
\hline 14&
\begin{tabular}{c}
$(x,A_1)$\\
$(y,A_1)$\\
$(z,A_1)$
\end{tabular}				&$(A_1^*)^{\oplus 4}$&
\begin{tabular}{l}
$P_1WL(x,y)$, $P_2WL(y,z)$, \\
$P_3WL(x,z)$, $P_4CT(x,y,z;-)$
\end{tabular}
\\
\hline 15&
$(x,A_4)$	&					
$\frac{1}{10}
\left[
\begin{array}{ccc}
3&1&-1\\
1&7&3\\
-1&3&7
\end{array}
\right]
$&
$P_1WL(x)$, $P_2WL(\eta_2,x)$, $P_3WL(\eta_3,x)$
\rule[-6mm]{0mm}{14mm}\\
\hline 16&
$(x,A_3)$			&$A_3^*\oplus A_1^*$&
\begin{tabular}{l}
$P_1WL(\eta_1,x)$, $P_2CT(-;x)$, \\
$P_3WL(\eta_3,x)$, $P_4WL(x)$
\end{tabular}
\\
\hline 17&
\begin{tabular}{c}
$(x,A_2)$\\
$(y,A_2)$
\end{tabular}					&$A_2^*\oplus \langle 1/6\rangle$&
\begin{tabular}{l}
$P_1CT(x,y;-)$, $P_2CT(x,y;-)$,
 $P_3WL(x,y)$
\end{tabular}
\\
\hline 18&
\begin{tabular}{c}
$(x,A_2)$\\
$(y,A_1)$
\end{tabular}					&		
$\frac{1}{6}
\left[
\begin{array}{cccc}
2&1&0&-1\\
1&5&3&1\\
0&3&6&3\\
-1&1&3&5
\end{array}
\right]
$&
\begin{tabular}{l}
$P_1WL(x,y)$, $P_2WL(\eta_2,x)$,\\ 
$P_3WL(\eta_3,y)$,
$P_4WL(\eta_4,x)$
\end{tabular}
\rule[-8mm]{0mm}{18mm}\\
\hline 19&
\begin{tabular}{c}
$(x,A_1)$\\
$(y,A_1)$
\end{tabular}
			&$D_4^*\oplus A_1^*$&
\begin{tabular}{l}		
$P_1WL(\eta_1,x)$,
$P_2WL(\eta_2,y)$,\\ 
$P_3CT(x,y;-)$, $P_4CT(-;-)$, $P_5WL(x,y)$
\end{tabular}			\\
\hline 20&
$(x,A_2)$			&$A_5^*$&
\begin{tabular}{l}
$P_1WL(\eta_1,x)$, $P_2CT(x;-)$,\\
$P_3BL$, $P_4CT(x;-)$, $P_5WL(\eta_5,x)$
\end{tabular}
\\
\hline 21&
$(x,A_1)$			&$D_6^*$&
\begin{tabular}{l}
$P_1WL(\eta_1,x)$, $P_2BL$,\\
$P_3CT(x;-)$, $P_4CT(-;-)$, \\
$s_{P_5}\cdot O=s_{P_5}\cdot \Theta_{x,1}=s_{P_5}\cdot \Theta_{\infty,1}=1$,\\
$s_{P_6}\cdot O=1$, $s_{P_6}\cdot \Theta_{x,1}=s_{P_6}\cdot \Theta_{\infty,1}=0$ 
\end{tabular}

\\
\hline
\end{longtable}

\subsection{Line-sections arising from weak-bitangents and bitangents and their description in $\MW(S_{\mcQ,z_o})$}\label{sec:line-sec}
Let $\mcQ$ be a singular quartic curve satisfying $(\dagger)$ and let $z_o$ be a general smooth point on $\mcQ$.

In order to prove Corollaries \ref{cor:ha3.3}, \ref{cor:ha3.4}, \ref{cor:ha3.5} and Theorem \ref{thm:A2+A1}, we consider bitangent lines and weak-bitangent lines written as $BL$ and $WL(\eta,x)$  as in Table \ref{table:geo-gene}, where $\eta$ is a smooth point of $\mcQ$ and $x$ is a singularity of $\mcQ$.
\begin{rem}
	If $\mcQ$ has a weak-bitangent line written as $WL(\eta,x)$ then $\mcQ$ must be an irreducible quartic curve which has double points only by the definition of weak-bitangent lines.
	In particular, the types of singular fibers of $S_{\mcQ,z_o}$ are only $\I_n$, $\III$ or $\IV$. 
\end{rem}
By Lemma \ref{lem:line-sec}, we obtain Proposition \ref{prop:w-bi-hieght}.
\begin{prop}\label{prop:w-bi-hieght}
	Let $\mcQ$ be an irreducible singular quartic curve with double points only and let $z_o$ be a general smooth point of $\mcQ$.
	Then, for $P\in E_{S_{\mcQ,z_o}}(\CC(t))$, the following conditions {\rm (i)} and {\rm (ii)} are equivalent:
	\begin{itemize}
		\item[{\rm (i)}]
			There exists some natural number $n_P$ such that $\langle P,P\rangle=3/2-n_P/(n_P+1)$  and the intersection number $s_P\cdot \Theta_{\infty,1}=1$.
		\item[{\rm (ii)}]
			$s_P$ is a line-section of a weak-bitangent line $L$ such that
			\begin{itemize}
				\item[{\rm (a)}] 
					$L$ passes through a singularity, $x_0$, of $\mcQ$, 
				\item[{\rm (b)}] 
					when the type of $x_0$ is $A_1$, $I_{x_0}(L,\mcQ)=2$ or $4$, i.e. $L$ is respectively the form $WL(\eta,x)$ or $WL(x)$ and
				\item[{\rm (c)}]
					otherwise, $I_{x_0}(L,\mcQ)=2$ i.e. $L$ is the form $WL(\eta,x)$.
			\end{itemize}
	\end{itemize}
\end{prop}
\proof
We assume {\rm (i)} and we choose $P\in E_{S_{\mcQ,z_o}}(\CC(t))$ satisfying (i).
Let $\sing(\mcQ)$ be the set of all singularities of $\mcQ$.  
We recall the formula of the height pairing as follows:
\[
	\langle P,P\rangle=2+2s_P\cdot O-\sum_{x\in\sing(\mcQ)\cup\{z_o\}}\mathrm{contr}_x(s_P,s_P)
\]
Observing $R_{\mcQ,z_o}$ in Table \ref{table:st-MW}, we see $\sum_{x\in\sing(\mcQ)\cup\{z_o\}}\mathrm{contr}_x(s_P,s_P)\leq 5/2$.
Hence we have 
\[
	\langle P,P\rangle=2+2s_P\cdot O-\sum_{x\in\sing(\mcQ)\cup\{z_o\}}\mathrm{contr}_x(s_P,s_P)\geq 2s_P\cdot O-1/2.
\]
By the condition {\rm (i)}, we have $s_P\cdot O=0$. 
Hence, $s_P$ is a line-section by Lemma \ref{lem:line-sec}. 
The following equality holds:
\[
	3/2-n_P/(n_P+1)=3/2-\sum_{x\in\sing(\mcQ)}\mathrm{contr}_x(s_P,s_P),
\]
so we have $\sum_{x\in\sing(\mcQ)}\mathrm{contr}_x(s_P,s_P)=n_P/(n_P+1)<1$.
For $x\in\sing(\mcQ)$, we must prove (ii) in the case when the type of $F_x$ is $\I_{m+1}$, $\III$ or $\IV$.
But we consider the case when the type of $F_x$ is $\I_{m+1} \ (m\geq1)$ since $\III$ and $\IV$ are special case of $\I_2$ and $\IV$, respectively
We write $F_x=\Theta_{x,0}+\sum^{m_x-1}_{i=1}a_{x,i}\Theta_{x,i}$.
Note that $m_x=m+1$.
If $s_P$ meets the $k$-th components $\Theta_{x,k} \ (0\leq k \leq m)$, then $\mathrm{contr}_x(s_P,s_P)=(m+1-k)k/(m+1)$.
In particular, 
\[
	\mathrm{contr}_x(s_P,s_P)
	\left\{
		\begin{array}{ll}
			=m/(m+1)& k=1,m\\
			=0& k=0\\
			>1& {\rm otherwise}
		\end{array}
	\right.
\]
Hence, by $\sum_{x\in\sing(\mcQ)}\mathrm{contr}_x(s_P,s_P)<1$, we have 
\[
	\mathrm{contr}_x(s_P,s_P)=0 \ {\rm or} \ m/(m+1)
\]
for all $x\in\sing(\mcQ)$.
Note that $m_1/(m_1+1)+m_2/(m_2+1)\geq 1$ where $m_1$ and $m_2$ are natural numbers. 
Therefore, there exists a unique $x_0\in \sing(\mcQ)$ such that
\begin{itemize}
	\item
	$\mathrm{contr}_{x_0}(s_P,s_P)=n_P/(n_P+1)$ and the type of $F_{x_0}$ is $\I_{n_P}$, and
	\item
	$\mathrm{contr}_{x}(s_P,s_P)=0$ for all $x\in \sing(\mcQ)\setminus \{x_0\}$. 
\end{itemize}
Therefore the type of $x_0$ is $A_{n_P}$ and we obtain $m_{x_0}=n_P+1$.

We assume that the type of $x_0$ is $A_1$.
Then we have $1/2=\mathrm{contr}_{x_0}(s_P,s_P)=s_P\cdot\Theta_{x_0,1}/2$.
By $s_P\cdot\Theta_{x_0,1}=1$, the intersection multiplicity of $\overline{f}_{\mcQ,z_o}(s_P)$ and a branch at $x_0$ of $\mcQ$ is one or more.
Considering that $I_{x_0}(\overline{f}_{\mcQ,z_o}(s_P),\mcQ)$ is even, it is 2 or 4.

When $I_{x_0}(\overline{f}_{\mcQ,z_o}(s_P),\mcQ)$ is 2, there exists a smooth $\eta$ on $\mcQ$ such that $I_{\eta}(\overline{f}_{\mcQ,z_o}(s_P),\mcQ)=2$.
Then we can write $L$ as $WL(\eta,x_0)$.

When $I_{x_0}(\overline{f}_{\mcQ,z_o}(s_P),\mcQ)$ is four, $\overline{f}_{\mcQ,z_o}(s_P)$ is tangent to a branch at $x_0$ of $\mcQ$ with multiplicity 3.
Then we can write $L$ as $WL(x_0)$

We next assume that $n_P>1$.
By consideration similar to the above case, $s_P\cdot\Theta_{x_0,1}=1$ or $s_P\cdot\Theta_{x_0,n_P}=s_P\cdot\Theta_{x_0,m_{x_0}-1}=1$.
For each branch at $x_0$ of $\mcQ$, the intersection multiplicities of $\overline{f}_{\mcQ,z_o}(s_P)$ are one.
Therefore we have $I_{x_0}(\overline{f}_{\mcQ,z_o}(s_P),\mcQ)=2$.
In the same way as the above case, the line $\overline{f}_{\mcQ,z_o}(s_P)$ is a weak-bitangent line as $WL(\eta,x_0)$, where $\eta$ is a smooth point of $\mcQ$.

Conversely, let $s_P$ and $x_0$ be a line-section and a singularity of $\mcQ$ satisfying the condition (ii).
We define the type of $x_0$ by $A_{n_P}$.
By Lemma \ref{lem:line-sec}, $s_{P}\cdot O=1$. 
By (b) and (c), we have $s_{P}\cdot \Theta_{x_0,1}=1$ or $s_{P}\cdot \Theta_{x_0,m_{x_0}-1}=1$.
We define $m_{x_0}-1$ by $n_P$.
Therefore, we obtain 
\[
	\mathrm{contr}_{x_0}(s_P,s_P)=(m_{x_0}-1)/m_{x_0}=n_P/(n_P+1). 
\]
Since $L$ is a weak-bitangent line satisfying (ii), we have $s_P\cdot\Theta_{x,0}=1$ for $x \in \sing(\mcQ)\setminus\{x_0\}$.
Hence we have $\langle P,P\rangle=3/2-n_P/(n_P+1)$.
\qed

Similarly we obtain the following proposition: 
\begin{prop}\label{prop:bi-height}
	Let $\mcQ$ be a plane quartic curve satisfying $(\dagger)$ and let $z_o$ be a general smooth point of $\mcQ$.
	Then, for $P\in E_{S_{\mcQ,z_o}}(\CC(t))$, the following conditions {\rm (i)} and {\rm (ii)} are equivalent:
	\begin{itemize}
		\item[{\rm (i)}]
			$\langle P,P\rangle=3/2$ and $s_P\cdot \Theta_{\infty,1}=1$.
		\item[{\rm (ii)}]
			$s_P$ is a line-section of a bitangent line or a $4$-fold tangent line.
	\end{itemize}
\end{prop}

By using Proposition \ref{prop:w-bi-hieght}, Remark \ref{rem:become-line} and Theorem \ref{lem:shio-form}, we obtain representations of weak-bitangent lines written as $WL(\eta,x)$ and bitangent lines. 
Note that a 4-fold tangent line (resp. a weak-bitangent line written as $WL(x)$ in the case when the type of $x$ is $A_1$) is the special case of a bitangent line (resp. a weak-bitangent line written as $WL(\eta,x)$) by Proposition \ref{prop:w-bi-hieght} and \ref{prop:bi-height}.
We explain some notations used in the next table.
\begin{itemize}
\item
	The columns labeled $\sing(\mcQ)$ and $E_{S_{\mcQ,z_o}}(\CC(t))$ are the same as those in Table \ref{table:st-MW} and $\eta$, $WL$ and $BL$ also are the same as those in  Table \ref{table:st-MW}.
\item
	Let $P_1,\ldots,P_n\in E_{S_{\mcQ,z_o}}(\CC(t))$ be generators as in Table \ref{table:geo-gene} and let $P$ be a rational point of $E_{S_{\mcQ,z_o}}(\CC(t))$ such that the corresponding section $s_P$ is a line-section of a bitangent line or a 4-fold tangent line.
	If the $P$ is given by $P_{s}=c_1P_1\dot{+}\cdots \dot{+}c_nP_n$, we write 
	$
		\left[
			\begin{array}{c}
				c_1\\
				\vdots\\
				c_n
			\end{array}
		\right]
	$
	for $PBL$ in the column labeled BL.
	Note that $s_{\dot{-}P}$ is also a line-section and $\overline{f}_{\mcQ,z_o}(s_{P})=\overline{f}_{\mcQ,z_o}(s_{\dot{-}P})$ but we do not write
$\left[
	\begin{array}{c}
		-c_1\\
		\vdots\\
		-c_n
	\end{array}
\right]$.

\item
	Let $P_1,\ldots,P_n\in E_{S_{\mcQ,z_o}}(\CC(t))$ be generators as Table \ref{table:geo-gene} and let $P_{\eta,x}$ be a rational point of $E_{S_{\mcQ,z_o}}(\CC(t))$ such that the corresponding section $s_{P_{\eta,x}}$ is a line-section of a weak-bitangent line written as $WL(\eta,x)$.
	If $P_{\eta,x}$ is given by $P_{\eta,x}=c^{\eta,x}_1P_1\dot{+}\cdots \dot{+}c^{\eta,x}_nP_n$, we write 
	$
	\left[
		\begin{array}{c}
			c^{\eta,x}_1\\
			\vdots\\
			c^{\eta,x}_n
		\end{array}
	\right]_{\eta,x}
	$
	for $P_{\eta,x}WL(\eta,x)$ in the column labeled WL.
	Note that $s_{\dot{-}P_{\eta,x}}$ is also a line-section and $\overline{f}_{\mcQ,z_o}(s_{P_{\eta,x}})=\overline{f}_{\mcQ,z_o}(s_{\dot{-}P_{\eta,x}})$ but we do not write
	$\left[
		\begin{array}{c}
			-c^{\eta,x}_1\\
			\vdots\\
			-c^{\eta,x}_n
			\end{array}
		\right]_{\eta,x}$.
	In the case when the type of $x$ is $A_1$, let $P_x=c^{\eta,x}_1P_1\dot{+}\cdots \dot{+}c^{\eta,x}_nP_n$ be a rational point corresponding to a line-section of a weak-bitangent line written as $WL(x)$.
	We also write 
	$
	\left[
		\begin{array}{c}
			c^{x}_1\\
			\vdots\\
			c^{x}_n
		\end{array}
	\right]_{\eta,x}
	$ for $P_{x}WL(x)$ in the column WL because $WL(x)$ is a special case of $WL(\eta,x)$.

\end{itemize}
\begin{center}
\footnotesize
\begin{longtable}{|c|c|c|c|}    
\caption{}
\label{table:geo-line}
\endhead
\hline
No. & $\sing(\mcQ)$ & BL & WL\\
\hline 1&
$(x,A_6)$		&
not exist&
$[3]_{\eta,x}$
\\
\hline	2&
$(x,E_6)$	&	
$[3]$&not exist
\\
\hline 3&
$(x,A_5)$	&	
$\left[\begin{array}{c}
0\\
3
\end{array}\right]$
&
$\left[\begin{array}{c}
1\\
1
\end{array}\right]_{\eta_1,x}$
$\left[\begin{array}{c}
-1\\
1
\end{array}\right]_{\eta_2,x}$
\rule[-4mm]{0mm}{9mm}\\
\hline 4&
$(x,D_5)$	&	
$\left[
\begin{array}{c}
1\\
\pm2
\end{array}
\right]$
&not exist
\\
\hline 5&
$(x,D_4)$	&		
\begin{tabular}{c}					
$\left[\begin{array}{c}
1\\
1\\
1\\
0\\
\end{array}\right]$
$\left[\begin{array}{c}
-1\\
1\\
1\\
0
\end{array}\right]$
$\left[\begin{array}{c}
1\\
-1\\
1\\
0
\end{array}\right]$
$\left[\begin{array}{c}
1\\
1\\
-1\\
0
\end{array}\right]$
\end{tabular}&not exist
\\
\hline 6&
\begin{tabular}{c}
$(x,A_4)$\\
$(y,A_2)$
\end{tabular}
	&	not exist&
$[3]_{\eta_1,x}$ $[5]_{\eta_2,y}$
\\
\hline 7&
\begin{tabular}{c}
$(x,A_4)$\\
$(y,A_1)$
\end{tabular}
		&

$\left[\begin{array}{c}
2\\
1
\end{array}\right]$&
$\left[\begin{array}{c}
-2\\
1
\end{array}\right]_{\eta_1,x}$
$\left[\begin{array}{c}
-1\\
2
\end{array}\right]_{\eta_2,y}$
\rule[-4mm]{0mm}{9mm}\\

\hline 8&

\begin{tabular}{c}
$(x,A_1)$\\
$(y,A_1)$\\
$(z,A_1)$\\
$(w,A_1)$
\end{tabular}
		&
			\begin{tabular}{c}
$\left[\begin{array}{c}
1\\
1\\
1\\
0\\
\end{array}\right]$
$\left[\begin{array}{c}
-1\\
1\\
1\\
0
\end{array}\right]$
$\left[\begin{array}{c}
1\\
-1\\
1\\
0
\end{array}\right]$
$\left[\begin{array}{c}
1\\
1\\
-1\\
0
\end{array}\right]$
\end{tabular}
&not exist

\\
\hline 9&
\begin{tabular}{c}
$(x,A_3)$\\
$(y,A_2)$
\end{tabular}
		&not exist&
\begin{tabular}{l}
$\left[\begin{array}{c}
0\\
3
\end{array}\right]_{\eta_1,x}$
$\left[\begin{array}{c}
1\\
2
\end{array}\right]_{\eta_2,y}$
$\left[\begin{array}{c}
-1\\
2
\end{array}\right]_{\eta_3,y}$
\end{tabular}
\rule[-4mm]{0mm}{9mm}\\
\hline 10&
\begin{tabular}{c}
$(x,A_3)$\\
$(y,A_1)$
\end{tabular}			&
$\left[\begin{array}{c}
0\\
\pm1\\
2
\end{array}\right]$&
\begin{tabular}{l}
$\left[\begin{array}{c}
1\\
0\\
1
\end{array}\right]_{\eta_1,x}$
$\left[\begin{array}{c}
-1\\
0\\
1
\end{array}\right]_{\eta_2,x}$
$\left[\begin{array}{c}
1\\
1\\
0
\end{array}\right]_{\eta_3,y}$\\
$\left[\begin{array}{c}
-1\\
1\\
0
\end{array}\right]_{\eta_4,y}$
\end{tabular}
\rule[-11mm]{0mm}{23mm}\\
\hline 11&
\begin{tabular}{c}
$(x,A_2)$\\
$(y,A_2)$\\
$(z,A_2)$
\end{tabular}
&
$\left[
\begin{array}{c}
3\\
0
\end{array}
\right]$&not exist
\\
\hline 12&
\begin{tabular}{c}
$(x,A_2)$\\
$(y,A_2)$\\
$(z,A_1)$
\end{tabular}
		&
$\left[
\begin{array}{c}
3\\
0
\end{array}
\right]$&

$\left[
\begin{array}{c}
-1\\
2
\end{array}
\right]_{\eta_1,x}$
$\left[
\begin{array}{c}
1\\
2
\end{array}
\right]_{\eta_2,y}$
\\

\hline 13&
\begin{tabular}{c}
$(x,A_2)$\\
$(y,A_1)$\\
$(z,A_1)$
\end{tabular}			&	
$\left[
\begin{array}{c}
1\\
1\\
1
\end{array}
\right]$
$\left[
\begin{array}{c}
1\\
-1\\
-1
\end{array}
\right]$
&
\begin{tabular}{l}
$\left[
\begin{array}{c}
1\\
-1\\
1
\end{array}
\right]_{\eta_1,x}$
$\left[
\begin{array}{c}
-1\\
-1\\
1
\end{array}
\right]_{\eta_2,x}$
$\left[
\begin{array}{c}
0\\
-1\\
2
\end{array}
\right]_{\eta_3,y}$\\
$\left[
\begin{array}{c}
0\\
-2\\
1
\end{array}
\right]_{\eta_4,z}$
\end{tabular}\rule[-11mm]{0mm}{23mm}\\
\hline 14&
\begin{tabular}{c}
$(x,A_1)$\\
$(y,A_1)$\\
$(z,A_1)$
\end{tabular}				&
			\begin{tabular}{l}
$\left[\begin{array}{c}
1\\
1\\
1\\
0\\
\end{array}\right]$
$\left[\begin{array}{c}
-1\\
1\\
1\\
0
\end{array}\right]$
$\left[\begin{array}{c}
1\\
-1\\
1\\
0
\end{array}\right]$
$\left[\begin{array}{c}
1\\
1\\
-1\\
0
\end{array}\right]$
\end{tabular}
&
\begin{tabular}{l}
$\left[\begin{array}{c}
1\\
0\\
0\\
1
\end{array}\right]_{\eta_1,z}$
$\left[\begin{array}{c}
1\\
0\\
0\\
-1
\end{array}\right]_{\eta_2,z}$

$\left[\begin{array}{c}
0\\
1\\
0\\
1
\end{array}\right]_{\eta_3,x}$\\
$\left[\begin{array}{c}
0\\
1\\
0\\
-1
\end{array}\right]_{\eta_4,x}$

$\left[\begin{array}{c}
0\\
0\\
1\\
1
\end{array}\right]_{\eta_5,y}$
$\left[\begin{array}{c}
0\\
0\\
1\\
-1
\end{array}\right]_{\eta_6,y}$
\end{tabular}
\rule[-15mm]{0mm}{31mm}\\
\hline 15&
$(x,A_4)$	&					
\begin{tabular}{l}
$\left[\begin{array}{c}
-1\\
-1\\
1
\end{array}\right]$
$\left[\begin{array}{c}
2\\
0\\
1
\end{array}\right]$
$\left[\begin{array}{c}
-2\\
1\\
0
\end{array}\right]$
\end{tabular}
&
$\left[\begin{array}{c}
0\\
1\\
0
\end{array}\right]_{\eta_1,x}$
$\left[\begin{array}{c}
0\\
0\\
1
\end{array}\right]_{\eta_2,x}$
$\left[\begin{array}{c}
1\\
-1\\
1
\end{array}\right]_{\eta_3,x}$
\rule[-6mm]{0mm}{14mm}\\
\hline 16&
$(x,A_3)$			&
\begin{tabular}{l}
$\left[\begin{array}{c}
0\\
1\\
0\\
\pm1
\end{array}\right]$
$\left[\begin{array}{c}
1\\
-1\\
1\\
\pm1
\end{array}\right]$
$\left[\begin{array}{c}
-1\\
0\\
1\\
\pm1
\end{array}\right]$
\end{tabular}
&
\begin{tabular}{l}
$\left[\begin{array}{c}
1\\
0\\
0\\
0
\end{array}\right]_{\eta_1,x}$
$\left[\begin{array}{c}
0\\
0\\
1\\
0
\end{array}\right]_{\eta_2,x}$
$\left[\begin{array}{c}
1\\
-1\\
0\\
0
\end{array}\right]_{\eta_3,x}$\\
$\left[\begin{array}{c}
0\\
-1\\
1\\
0
\end{array}\right]_{\eta_4,x}$
\end{tabular}
\rule[-15mm]{0mm}{31mm}\\
\hline 17&
\begin{tabular}{c}
$(x,A_2)$\\
$(y,A_2)$
\end{tabular}					&
$\left[\begin{array}{c}
0\\
0\\
3
\end{array}\right]$&
\begin{tabular}{c}
$\left[\begin{array}{c}
1\\
0\\
1
\end{array}\right]_{\eta_1,x}$
$\left[\begin{array}{c}
0\\
1\\
-1
\end{array}\right]_{\eta_2,x}$
$\left[\begin{array}{c}
1\\
-1\\
-1
\end{array}\right]_{\eta_3,x}$\\
$\left[\begin{array}{c}
1\\
0\\
-1
\end{array}\right]_{\eta_4,y}$
$\left[\begin{array}{c}
0\\
1\\
1
\end{array}\right]_{\eta_5,y}$
$\left[\begin{array}{c}
1\\
-1\\
1
\end{array}\right]_{\eta_6,y}$
\end{tabular}
\rule[-11mm]{0mm}{23mm}\\
\hline 18&
\begin{tabular}{c}
$(x,A_2)$\\
$(y,A_1)$
\end{tabular}					&		
\begin{tabular}{l}
$\left[\begin{array}{c}
2\\
0\\
0\\
1
\end{array}\right]$
$\left[\begin{array}{c}
-1\\
0\\
-1\\
1
\end{array}\right]$
$\left[\begin{array}{c}
1\\
1\\
-1\\
0
\end{array}\right]$
$\left[\begin{array}{c}
-2\\
1\\
0\\
0
\end{array}\right]$
\end{tabular}
&
\begin{tabular}{l}
$\left[\begin{array}{c}
0\\
1\\
0\\
0
\end{array}\right]_{\eta_1,x}$
$\left[\begin{array}{c}
0\\
0\\
0\\
1
\end{array}\right]_{\eta_2,x}$
$\left[\begin{array}{c}
1\\
0\\
-1\\
1
\end{array}\right]_{\eta_3,x}$\\
$\left[\begin{array}{c}
-1\\
1\\
-1\\
0
\end{array}\right]_{\eta_4,x}$
$\left[\begin{array}{c}
0\\
0\\
1\\
0
\end{array}\right]_{\eta_5,y}$
$\left[\begin{array}{c}
0\\
1\\
-1\\
1
\end{array}\right]_{\eta_6,y}$\\
$\left[\begin{array}{c}
1\\
-1\\
0\\
1
\end{array}\right]_{\eta_7,y}$

\end{tabular}
\rule[-22mm]{0mm}{44mm}\\
\hline 19&
\begin{tabular}{c}
$(x,A_1)$\\
$(y,A_1)$
\end{tabular}&
\begin{tabular}{l}
$\left[\begin{array}{c}
0\\
0\\
1\\
0\\
1
\end{array}\right]$
$\left[\begin{array}{c}
0\\
0\\
1\\
-1\\
1
\end{array}\right]$
$\left[\begin{array}{c}
0\\
0\\
-1\\
0\\
1
\end{array}\right]$
$\left[\begin{array}{c}
0\\
0\\
-1\\
1\\
1
\end{array}\right]$\\
$\left[\begin{array}{c}
-1\\
1\\
0\\
0\\
1
\end{array}\right]$
$\left[\begin{array}{c}
1\\
1\\
0\\
-1\\
1
\end{array}\right]$
$\left[\begin{array}{c}
1\\
-1\\
0\\
0\\
1
\end{array}\right]$
$\left[\begin{array}{c}
-1\\
-1\\
0\\
1\\
1
\end{array}\right]$
\end{tabular}&
\begin{tabular}{l}
$\left[\begin{array}{c}
1\\
0\\
0\\
0\\
0
\end{array}\right]_{\eta_1,x}$
$\left[\begin{array}{c}
0\\
1\\
1\\
-1\\
0
\end{array}\right]_{\eta_2,x}$
$\left[\begin{array}{c}
0\\
-1\\
1\\
0\\
0
\end{array}\right]_{\eta_3,x}$\\
$\left[\begin{array}{c}
-1\\
0\\
0\\
1\\
0
\end{array}\right]_{\eta_4,x}$

$\left[\begin{array}{c}
0\\
1\\
0\\
0\\
0
\end{array}\right]_{\eta_5,y}$
$\left[\begin{array}{c}
-1\\
0\\
1\\
0\\
0
\end{array}\right]_{\eta_6,y}$\\
$\left[\begin{array}{c}
1\\
0\\
1\\
-1\\
0
\end{array}\right]_{\eta_7,y}$
$\left[\begin{array}{c}
0\\
1\\
0\\
-1\\
0
\end{array}\right]_{\eta_8,y}$
\end{tabular}
\rule[-27mm]{0mm}{45mm}\\
\hline 20&
$(x,A_2)$			&
\begin{tabular}{l}
$\left[\begin{array}{c}
1\\
-1\\
1\\
-1\\
1
\end{array}\right]$
$\left[\begin{array}{c}
1\\
-1\\
1\\
0\\
-1
\end{array}\right]$
$\left[\begin{array}{c}
1\\
0\\
-1\\
0\\
1
\end{array}\right]$
$\left[\begin{array}{c}
1\\
0\\
-1\\
1\\
-1
\end{array}\right]$\\
$\left[\begin{array}{c}
1\\
-1\\
0\\
1\\
0
\end{array}\right]$
$\left[\begin{array}{c}
1\\
0\\
0\\
-1\\
0
\end{array}\right]$
$\left[\begin{array}{c}
0\\
0\\
1\\
0\\
0
\end{array}\right]$
$\left[\begin{array}{c}
0\\
-1\\
1\\
-1\\
0
\end{array}\right]$\\
$\left[\begin{array}{c}
0\\
-1\\
0\\
0\\
1
\end{array}\right]$
$\left[\begin{array}{c}
0\\
1\\
0\\
-1\\
1
\end{array}\right]$
\end{tabular}&

\begin{tabular}{c}
$\left[\begin{array}{c}
1\\
0\\
0\\
0\\
0
\end{array}\right]_{\eta_1,x}$
$\left[\begin{array}{c}
1\\
-1\\
0\\
0\\
0
\end{array}\right]_{\eta_2,x}$
$\left[\begin{array}{c}
0\\
-1\\
1\\
0\\
0
\end{array}\right]_{\eta_3,x}$\\
$\left[\begin{array}{c}
0\\
0\\
1\\
-1\\
0
\end{array}\right]_{\eta_4,x}$
$\left[\begin{array}{c}
0\\
0\\
0\\
0\\
1
\end{array}\right]_{\eta_5,x}$
$\left[\begin{array}{c}
0\\
0\\
0\\
-1\\
1
\end{array}\right]_{\eta_6,x}$\\
\end{tabular}
\\
\hline 21&

$(x,A_1)$			&
\begin{tabular}{l}
$\left[\begin{array}{c}
1\\
1\\
0\\
-1\\
1\\
-1
\end{array}\right]$
$\left[\begin{array}{c}
1\\
1\\
0\\
0\\
-1\\
0
\end{array}\right]$
$\left[\begin{array}{c}
1\\
-1\\
0\\
-1\\
1\\
0
\end{array}\right]$
$\left[\begin{array}{c}
1\\
-1\\
0\\
0\\
-1\\
1
\end{array}\right]$\\
$\left[\begin{array}{c}
1\\
0\\
-1\\
-1\\
0\\
1
\end{array}\right]$
$\left[\begin{array}{c}
1\\
0\\
1\\
0\\
0\\
-1
\end{array}\right]$
$\left[\begin{array}{c}
1\\
0\\
-1\\
0\\
0\\
0
\end{array}\right]$
$\left[\begin{array}{c}
1\\
0\\
1\\
-1\\
0\\
0
\end{array}\right]$\\
$\left[\begin{array}{c}
0\\
1\\
0\\
1\\
0\\
-1
\end{array}\right]$
$\left[\begin{array}{c}
0\\
1\\
0\\
-1\\
0\\
0
\end{array}\right]$
$\left[\begin{array}{c}
0\\
1\\
0\\
0\\
0\\
0
\end{array}\right]$
$\left[\begin{array}{c}
0\\
1\\
0\\
0\\
0\\
-1
\end{array}\right]$\\
$\left[\begin{array}{c}
0\\
0\\
1\\
0\\
1\\
-1
\end{array}\right]$
$\left[\begin{array}{c}
0\\
0\\
1\\
-1\\
1\\
-1
\end{array}\right]$
$\left[\begin{array}{c}
0\\
0\\
-1\\
0\\
1\\
0
\end{array}\right]$
$\left[\begin{array}{c}
0\\
0\\
-1\\
-1\\
1\\
0
\end{array}\right]$
\end{tabular}&

\begin{tabular}{c}
$\left[\begin{array}{c}
1\\
0\\
0\\
0\\
0\\
0
\end{array}\right]_{\eta_1,x}$
$\left[\begin{array}{c}
1\\
0\\
0\\
-1\\
0\\
0
\end{array}\right]_{\eta_2,x}$
$\left[\begin{array}{c}
0\\
1\\
1\\
0\\
0\\
-1
\end{array}\right]_{\eta_3,x}$\\
$\left[\begin{array}{c}
0\\
1\\
-1\\
0\\
0\\
0
\end{array}\right]_{\eta_4,x}$
$\left[\begin{array}{c}
0\\
0\\
0\\
-1\\
1\\
0
\end{array}\right]_{\eta_5,x}$
$\left[\begin{array}{c}
0\\
0\\
0\\
0\\
1\\
-1
\end{array}\right]_{\eta_6,x}$
\end{tabular}
\\
\hline

\end{longtable}
\end{center}

\section{The Mumford representation of a semi-reduced divisor}\label{sec:mum-rep}
\subsection{Mumford representations}
Let $C$ be a hyperelliptic curve of genus $g$ defined over a field $K \ (\mathrm{char}(K)\neq2)$ given by an affine equation:
\[
    y^2=f(x)=x^{2g+1}+c_1x^{2g}+\cdots +c_{2g+1} \quad c_i\in K \ (i=1,\ldots ,2g+1).
\]
Its point at infinite is denoted by $O$ and the hyperelliptic involution by $\iota:C\rightarrow C$.
\begin{defin}\label{def:mumford}
    For a divisor $\fd=\sum_{P\in C}n_PP\in \Div(C)$, we assume that $O \notin \Supp(\fd):=\{ P\in C \mid n_P\neq 0 \}$ and $\fd$ is effective.
    We call $\fd$ a {\it semi-reduced divisor} if it follows that
    \begin{itemize}
        \item if $P\in\Supp(\fd)$ and $P\neq\iota(P)$, then $\iota(P)\notin \Supp(\fd)$ and
        \item if $P\in\Supp(\fd)$ and $P=\iota(P)$, then $n_P=1$.
    \end{itemize}
\end{defin}
We denote the coordinate ring $\Kbar[x,y]/\langle y^2-f \rangle$ by $\Kbar[C]$ and the image of $g\in\Kbar[x,y]$ in $\Kbar[C]$ by $[g]$. For $P\in C$, we write $\mcO_P$ for the local ring at $P$ and $\ord_P$ for the discrete valuation at $P$.
Then we have
\begin{prop}{\rm (\cite[Proposition 2.1]{Taka-Toku})}
    Let $>_p$ be the pure lexicographical order with $y>_px$ in $\Kbar[x,y]$.  For a semi-reduced divisor $\fd=\sum_{P\in C}n_PP$, we put
    \[
        \widetilde{I(\fd)}:=\{ g\in \Kbar[C] \mid \ord_P([g])\geq n_P \ \forall P\in \Supp(\fd)  \}.
    \]
    Then the reduced Gr\"{o}bner basis of $\widetilde{I(\fd)}$ with respect to $>_p$ is of the form $\{ a(x) , y-b(x) \}$ and $b^2-f\in \langle a \rangle$, where $a$ and $b\in\Kbar[x]$. 
\end{prop}
For a proof, see \cite[Proposition 2.1]{Taka-Toku}.

\begin{defin}
    Let $\fd$ be a semi-reduced divisor on $C$ and let $\{ a(x) , y-b(x) \}$ be the reduced Gr\"{o}bner basis for $\widetilde{I(\fd)}$ with respect to $>_p$. 
    Then we call the pair $(a,b)$ the {\it Mumford representation} of $\fd$.
\end{defin}

The next lemma is a characterization of Mumford representations. 
\begin{lem}{\rm (\cite[Lemma 10.3.5]{galb})}\label{lem:chara-Mum}
    Let $\fd=\sum_{P\in C}n_PP$ be a semi-reduced divisor and we put $P=(x_P,y_P)$.
    The pair $(a,b)\in (\Kbar[x])^2$ is the Mumford representation of $\fd$ if and only if $(a,b)$ satisfies
    \begin{equation*}
        {\rm(i)} \ a(x)\!=\!\mathop{\prod}\limits_{P\in \Supp(\fd)}(x-x_P)^{n_P}, \quad
        {\rm(ii)} \ \deg \,b(x) < \deg \,a(x), \quad 
        {\rm(iii)} \ a\mid b^2-f.
    \end{equation*}
\end{lem}
\subsection{Semi-reduced divisors of degree 3 on elliptic curves}
From now on, we assume that $K$ is a field of  characteristics $\neq2$. Let $E$ be an elliptic curve defined over $K$ given by the Weierstrass form
\[
    y^2=f(x)=x^3+c_1x^2+c_2x+c_3 \quad c_i\in K  \ (i=1,2,3).
\]
\begin{rem}
    {\rm (i)}
    If a semi-reduced divisor $\fd$ is defined over $K$, then the Mumford representation $(a,b)$ of $\fd$ belongs to ${(K[x])}^2$. 
    {\rm (ii)} 
    For a semi-reduced divisor  $\sum_in_iP_i$ defined over $K$, the points $P_i$'s are not always elements in $E(K)$.     
\end{rem}

Let $\fd=P_1+P_2+P_3$ be a semi-reduced divisor of degree 3 defined over $K$.
We write  $P_1\dot{+}P_2\dot{+}P_3$ for $P_{\fd}$, where $\dot{+}$ is the addition on $E$ and the infinity point is the identity with respect to $\dot{+}$. 
\begin{lem}{\rm (\cite[Lemma 5.6]{Ban-Kaw-Mas-Tok})}\label{lem:mum-tri1}
    Assume that $P_{\fd}\neq O$ and let $(a,b)$ be the Mumford representation of $\fd$. Then we have
    \begin{itemize}
        \item[{\rm(i)}] $P_{\fd}\neq P_i \ (i=1,2,3)$.
        \item[{\rm(ii)}] $\deg b=2$.
    \end{itemize}
\end{lem}
\begin{lem}{\rm (\cite[Lemma 5.7]{Ban-Kaw-Mas-Tok})}
    \label{lem:mum-tri2}
    We keep the notation as before.
    Assume that $\fd$ is defined over $K$. 
    Put $P_{\fd}:=(x_{\fd},y_{\fd})$. 
    Then we have the following:
    \begin{itemize}
        \item[{\rm(i)}] The point $P_{\fd}$ is a $K$-rational point of $E$.
        \item[{\rm(ii)}] The two polynomials $a$ and $b$ belong to $K[x]$ and, $b$ is of the form 
        \[
            b_0(x-x_{\fd})(x-b_1)-y_{\fd} \ (b_0,b_1\in K).
        \]
    \end{itemize}
\end{lem}

\section{The Proof of our main results}\label{sec:appli}

Let $\mcQ$ be a reduced quartic curve satisfying $(\dagger)$ and let $z_o$ be a general smooth point of $\mcQ$.
Let $s_1$, $s_2$ and $s_3$ be line-sections of $S_{\mcQ,z_o}$.
By using the Mumford representation of the semi-reduced divisor $P_{s_1}+P_{s_2}+P_{s_3}$,
we consider the geometry of its weak-bitangent lines, conics and $\mcQ$.
For this purpose, we need to give a Weierstrass equation of $E_{S_{\mcQ,z_o}}$ over $\CC(t)$.

In Section \ref{sec:setting}, we explain a Weierstrass equation which we need to consider.
Choose homogeneous coordinates $[T,X,Z]$ of $\PP^2$ such that $z_o=[0,1,0]$ and the tangent line at $z_o$ is given by $Z=0$ and 
$\mcQ$ is given by
\[
    F(T,X,Z)=X^3Z+A_2(T,Z)X^2+A_3(T,Z)X+A_4(T,Z),
\]
where $A_i(T,Z) \ i=2,3,4$ are homogeneous polynomials of degree $i$.
We denote affine coordinates by $(t,x)=(T/Z,X/Z)$.
In Section \ref{sec:proof_main1}, \ref{sec:proof_cors} and \ref{sec:proof_main2} we will consider the Weierstrass equation
\[
    y^2=F(t,x,1)
\]

\subsection{Settings}\label{sec:setting}
In order to prove Theorem \ref{thm:main}, we will prove the next lemmas.
Let $[T,X,Z]$  be homogeneous coordinates of $\PP^2$ and let $(t,x)=(T/Z,X/Z)$ be affine coordinates.
\begin{lem}\label{prop:set-qua}
    Let $\mcQ$ be a reduced quartic curve that is not four lines. 
    Then there is a coordinate system $[T,X,Z]$ of $\PP^2$ such that $\mcQ$ is given by
    \[
        X^3Z+A_{\mcQ,2}(T,Z)X^2+A_{\mcQ,3}(T,Z)X+A_{\mcQ,4}(T,Z)
    \]
    where $A_{\mcQ,d}$ is a binary form of degree $d$ in $T$ and $Z$ such that
    \[
        \deg\,A_{\mcQ,2}(t,1)=2, \ \deg\,A_{\mcQ,3}(t,1)=3, \ {\it and} \ \deg\,A_{\mcQ,4}(t,1)\leq 3.
    \]
\end{lem}
{\it Proof.} Let $z_o$ be a general smooth point on $\mcQ$ and let $L_{z_o}$ be the tangent line to $\mcQ$ at $z_o$.
We take a point $p\in \{p,p'\}:=(\mcQ\cap L_{z_o})\setminus \{ z_o \}$ and $q\in \mcQ \setminus L_{z_o}$.
Note that we may assume that $I_{p}(\mcQ,L_{z_o})=I_{p'}(\mcQ,L_{z_o})=1$ and $I_{z_o}(\mcQ,L_{z_o})=2$.
We can choose a coordinate system $[T,X,Z]$ such that
\begin{itemize}
    \item the point $q$ has coordinates $[0,0,1]$ and the point $p$ has coordinates $[1,0,0]$ and
    \item the smooth point $z_o$ has coordinates $[0,1,0]$ and the tangent line $L_{z_o}$ is defined by $Z=0$.
\end{itemize}
We have the desired coordinate system.
\qed

By comparing degree of polynomials, we have the next lemma. 
\begin{lem}\label{lem:degree}
Let $f(t,x)$ be a polynomial in $\CC[t,x]$ given by 
    \[
        f=x^3+a_2(t)x^2+a_3(t)x+a_4(t) \ (a_i\in\CC[t], \ i=2,3,4)
    \]
    and
    \[
        \deg\,a_2=2, \ \deg\,a_3=3 \ and \ \deg\,a_4\leq3.
    \]
    Let $x_0$ and $y_0$ be polynomials in $\CC[t]$ such that ${y_0}^2=f(t,x_0)$.
    For $(r,s)\in\CC(t)^{\times}\times\CC(t)$, we put $v_{r,s}=r(x-x_0)(x-s)-y_0$. 
    Note that it follows that ${v_{r,s}}^2-f=r^2(x-x_0)u$, where $u=x^3+c_1(t)x^2+c_2(t)x^2+c_3(t), \ (c_i\in\CC(t), \ i=1,2,3)$.

    If $u$ is a polynomial of total degree of $3$, then $r\in\CC^{\times}$ and $s\in\CC[t]$ of degree $\leq1$.
\end{lem}
\subsection{Proof of Theorem\,\ref{thm:main}}\label{sec:proof_main1}

From Lemma \ref{prop:set-qua}, it follows that there is a coordinate system $[T,X,Z]$ of $\PP^2$ such that $\mcQ$ is given by 
\[
    F_{\mcQ}(T,X,Z)=X^3Z+A_{\mcQ,2}(T,Z)X^2+A_{\mcQ,3}(T,Z)X+A_{\mcQ,4}(T,Z),
\]
where $(t,x)=(T/Z,X/Z)$ be non-homogeneous coordinates and $A_{\mcQ,d}$ are binary forms of degree $d$ in $T$ and $Z$ such that
\begin{align*}
    \deg\,A_{\mcQ,2}(t,1)=2, \ \deg\,A_{\mcQ,3}(t,1)=3 \ {\rm and } \ \deg\,A_{\mcQ,4}(t,1)\leq 3.
\end{align*}

We put $L_i=\overline{f}_{\mcQ,z_o}(s_i) \ (i=1,2,3)$ and $L_4=\overline{f}_{\mcQ,z_o}(s_{P_4})$.
The elliptic curve $E$ given by $y^2=F_{\mcQ}(t,x,1)$ corresponds to the generic fiber of $\varphi_{\mcQ,z_o}$.
The lines $L_i$ are given by the affine equations $x-x_i(t)=0 \ (i=1,2,3,4)$, where $x_1(t),x_2(t),x_3(t)$ and $x_4(t)$ are polynomials of degree $\leq$ 1 and we denote $P_1+P_2+P_3\in\Div(E)$ by $\fd$.
The lines $L_1$, $L_2$ and $L_3$ are three distinct lines, so $\iota(P)\notin \Supp(\fd)$ for $P\in\Supp(\fd)$.
Hence, the divisor $\fd$ is a semi-reduced divisor and we denote the Mumford representation of $\fd$ by $(a,b) \ (a,b\in\CC(t)[x])$.
From Lemma \ref{lem:chara-Mum} and Lemma \ref{lem:mum-tri2}, we have
\begin{align}\label{eq:tri}
    \begin{aligned}
        a&=(x-x_1)(x-x_2)(x-x_3),\\
        b&=b_0(x-x_{4})(x-b_1)-y_{4} \ {\rm and} \\
        b^2&-F_{\mcQ}(t,x,1)={b_0}^2(x-x_4)a,\\
    \end{aligned}
\end{align}
where $b_0\in \CC(t)^{\times}$ and $b_1\in\CC(t)$.
Now $a$ is a polynomial of total degree 3.
From Lemma \ref{lem:degree}, it follows that $b_0\in\CC^{\times}$ and $b_1\in\CC[t]$ such that $\deg\, b_1\leq 1$.
Therefore the affine equation $b(t,x)=0$ gives rise to a plane curve $C$ of degree $2$ i.e. $C$ is a conic.
From the third equation of (\ref{eq:tri}), $C$ is the desired conic.
\qed

\subsection{Proof of Corollaries \ref{cor:ha3.3}, \ref{cor:ha3.4} and \ref{cor:ha3.5}}\label{sec:proof_cors}

By Theorem \ref{thm:main} and Table \ref{table:geo-line}, we obtain Corollaries \ref{cor:ha3.3}, \ref{cor:ha3.4} and \ref{cor:ha3.5}.
We prove Corollary \ref{cor:ha3.3}, because we can similarly prove the other corollaries.
We keep the notation as before.

{\it Proof of Corollary \ref{cor:ha3.3}.}
{
Let $z_o$ be a general smooth point on $C_1+C_2$ and let $L_1,L_2,L_3$ and $L_4$ be the four distinct bitangent lines of $C_1+C_2$. We denote generators of $E_{S_{C_1+C_2,z_o}}(\CC(t))$ by $P_1,P_2,P_3$ and $P_4$ as Table \ref{table:geo-gene}.
By Table \ref{table:geo-line}, we take $Q_1,Q_2,Q_3$ and $Q_4$ as follows:
\begin{itemize}
    \item $Q_1:=\dot{-}P_1\dot{+}P_2\dot{+}P_3$,
    \item $Q_2:=P_1\dot{-}P_2\dot{+}P_3$,
    \item $Q_3:=P_1\dot{+}P_2\dot{-}P_3$,
    \item $Q_4:=P_1\dot{+}P_2\dot{+}P_3$ and
    \item $\overline{f}_{C_1+C_2,z_o}(s_{Q_i})=L_i \ (i=1,2,3,4)$.
\end{itemize}
Then $Q_4=Q_1\dot{+}Q_2\dot{+}Q_3$ holds.
From Theorem \ref{thm:main}, the eight points of $C_1+C_2$ of contact with the four lines lie on a conic $C$.
We shall prove take four distinct $p_1,p_2,p_3$ and $p_4$ of $(\cup_{i=1}^4L_i)\cap\mcQ$ such that the all three points of $p_1,p_2,p_3$ and $p_4$ do not lie on $C$.
Hence, $C$ is smooth.
\qed

\subsection{Proof of Theorem \ref{thm:A2+A1}}\label{sec:proof_main2}
Let $P_1,P_2,P_3$ and $P_4$ be generators of $E_{S_{\mcQ,z_o}}(\CC(t))$ as Table \ref{table:geo-gene}. 
From Table \ref{table:geo-line}, there exists seven lines $L_l$ and $M_m \ (l=1,2,3,4, \ m=1,2,3)$ satisfying as follows:
\begin{itemize}
    \item
    the lines $L_l$ is a weak-bitangent line passing through $x$ and tangent to $\mcQ$ at a smooth point for $l=1,2,3,4$ and
    \item
    the lines $M_m$ is a weak-bitangent line passing through $y$ and tangent to $\mcQ$ at a smooth point for $m=1,2,3$ and
    \item
    $\overline{f}_{\mcQ,z_o}(s_{Q_l})=L_l$ for $l=1,2,3,4$ and
    \item
    $\overline{f}_{\mcQ,z_o}(s_{R_m})=M_m$ for $m=1,2,3$.
\end{itemize}
where 
$Q_1=P_2$, $Q_2=P_4$, $Q_3=P_1\dot{-}P_3\dot{+}P_4$, $Q_4=\dot{-}P_1\dot{+}P_2\dot{-}P_3$, 
$R_1=P_3$, 
$R_2=P_2\dot{-}P_3\dot{+}P_4$ and $R_3=P_1\dot{-}P_2\dot{+}P_4$.
In order to explain clearly, we use similar notations of Table \ref{table:geo-line}.
Namely, we write
\[
    Q_1=\left[\begin{array}{c}
        0\\
        1\\
        0\\
        0
    \end{array}\right],
    Q_2=\left[\begin{array}{c}
        0\\
        0\\
        0\\
        1
    \end{array}\right],
    Q_3=\left[\begin{array}{c}
        1\\
        0\\
        -1\\
        1
    \end{array}\right],
    Q_4=\left[\begin{array}{c}
        -1\\
        1\\
        -1\\
        0
    \end{array}\right]
\]
\[
    R_1=\left[\begin{array}{c}
        0\\
        0\\
        1\\
        0
    \end{array}\right],
    R_2=\left[\begin{array}{c}
        0\\
        1\\
        -1\\
        1
    \end{array}\right],
    R_3=\left[\begin{array}{c}
        1\\
        -1\\
        0\\
        1
    \end{array}\right].
\]
We prove the existence of $M_{a_{ij}}$ and $M_{b_{ij}}$ $(1\leq a_{ij}<b_{ij}\leq 4)$ for $L_i$ and $L_j \ (1\leq i<j \leq 4)$.
We consider the case when $i=1$ and $j=2$.
Then we obtain $R_2=Q_1\dot{+}Q_2\dot{-}R_1$.
From Theorem \ref{thm:main}, the six points of $\mcQ$ of intersections with the four lines lie on a conic $C_{12}$.
Note that $C_{12}$ is smooth.
For other pair $(i,j)$, we similarly find the existence of $M_{a_{ij}}$, $M_{b_{ij}}$ and $C_{ij}$ satisfying $(\ast)$.

We will prove the uniqueness.
For $L_i$ and $L_j$, let $M_{a_{ij}}, M_{b_{ij}}$ and $C_{ij}$ be two lines and a smooth conic satisfying $(\ast)$ and put $\{ x, p_1, p_2\}=\mcQ \cap (L_i+L_j)$,  $\{ y, q_1\}=\mcQ \cap M_{a_{ij}}$ and $\{ y, q_2\}=\mcQ \cap M_{b_{ij}}$.
Note that the divisor on $C_{ij}$ cut out by $\mcQ$ is $C_{ij}|_{\mcQ}=2x+2y+p_1+p_2+q_1+q_2$.

Assume that for a line $M_{a_{ij}}$ and $M_{c_{ij}}$ $(c_{ij}\in \{1,2,3\} \setminus \{ a_{ij}\})$, there exists a smooth conic $C'_{ij}$ such that $(L_i+L_j+M_{a_{ij}}+M_{b_{ij}})\cap \mcQ \subset C'_{ij}$.
We will prove that $M_{b_{ij}}=M_{c_{ij}}$.

Put $\{ y,q_3\}=\mcQ\cap M_{c_{ij}}$.
Note that the divisor on $C'_{ij}$ cut out by $\mcQ$ is $C'_{ij}|_{\mcQ}=2x+2y+p_1+p_2+q_1+q_3$.

We have $\{x,y,p_1,p_2,q_1\}=\mcQ \cap (L_i+L_j+M_{a_{ij}})\subset C_{ij}, C'_{ij} \ i.e. \ C_{ij}=C'_{ij}$.
We obtain $2x+2y+p_1+p_2+q_1+q_2=C_{ij}|_{\mcQ}=C'_{ij}|_{\mcQ}=2x+2y+p_1+p_2+q_1+q_3$.
Therefore, $q_2=q_3$.
We have $M_{b_{ij}}=M_{c_{ij}}$.
\qed
\begin{rem}
    We proved Theorem \ref{thm:A2+A1} in the case when the two branch at $y$ of $\mcQ$ meets $L_i$ at $y$ with multiplicities 1, respectively.
    We also consider the case when a branch at $y$ of $\mcQ$ meets $L_i$ at $y$ with multiplicity 3.
    In this case, we can prove a similar statement.
\end{rem}

\bibliographystyle{abbrv}
\bibliography{ref}

\begin{thebibliography}{10}

\bibitem{Ban-Kaw-Mas-Tok}
S.~Bannai, N.~Kawana, R.~Masuya, and H.~Tokunaga.
\newblock Trisections on certain rational elliptic surfaces and families of
  {Z}ariski pairs degenerating to the same conic-line arrangement, 2021.
\newblock arXiv:2103.07639.

\bibitem{Ban-Toku15}
S.~Bannai and H.-o. Tokunaga.
\newblock Geometry of bisections of elliptic surfaces and {Z}ariski {$N$}-plets
  for conic arrangements.
\newblock {\em Geom. Dedicata}, 178:219--237, 2015.

\bibitem{Ban-Toku17}
S.~Bannai and H.-o. Tokunaga.
\newblock Geometry of bisections of elliptic surfaces and {Z}ariski {$N$}-plets
  {II}.
\newblock {\em Topology Appl.}, 231:10--25, 2017.

\bibitem{Ban-Toku21}
S.~Bannai and H.-o. Tokunaga.
\newblock Elliptic surfaces of rank one and the topology of cubic-line
  arrangements.
\newblock {\em J. Number Theory}, 221:174--189, 2021.

\bibitem{BHPV}
W.~P. Barth, K.~Hulek, C.~A.~M. Peters, and A.~Van~de Ven.
\newblock {\em Compact complex surfaces}, volume~4 of {\em Ergebnisse der
  Mathematik und ihrer Grenzgebiete. 3. Folge. A Series of Modern Surveys in
  Mathematics}.
\newblock Springer-Verlag, Berlin, second edition, 2004.

\bibitem{galb}
S.~D. Galbraith.
\newblock {\em Mathematics of public key cryptography}.
\newblock Cambridge University Press, Cambridge, 2012.

\bibitem{Harris}
J.~Harris.
\newblock \textup{Theta-characteristics on algebraic curves}.
\newblock {\em \textit{Trans. Amer. Math. Soc.}}, 271(2):611--638, 1982.

\bibitem{Hori}
E.~Horikawa.
\newblock On deformations of quintic surfaces.
\newblock {\em Invent. Math.}, 31(1):43--85, 1975.

\bibitem{Koda}
K.~Kodaira.
\newblock On compact analytic surfaces. {II}, {III}.
\newblock {\em Ann. of Math. (2) 77 (1963), 563--626; ibid.}, 78:1--40, 1963.

\bibitem{Mira}
R.~Miranda.
\newblock {\em The basic theory of elliptic surfaces}.
\newblock Dottorato di Ricerca in Matematica. [Doctorate in Mathematical
  Research]. ETS Editrice, Pisa, 1989.

\bibitem{Ogu-Shio}
K.~Oguiso and T.~Shioda.
\newblock \textup{The {M}ordell-{W}eil lattice of a rational elliptic surface}.
\newblock {\em \textit{Comment. Math. Univ. St. Paul.}}, 40(1):83--99, 1991.

\bibitem{Shio}
T.~Shioda.
\newblock \textup{On the {M}ordell-{W}eil lattices}.
\newblock {\em \textit{Comment. Math. Univ. St. Paul.}}, 39(2):211--240, 1990.

\bibitem{Shio93}
T.~Shioda.
\newblock Plane quartics and {M}ordell-{W}eil lattices of type {$E_7$}.
\newblock {\em Comment. Math. Univ. St. Paul.}, 42(1):61--79, 1993.

\bibitem{Taka-Toku}
A.~Takahashi and H.-o. Tokunaga.
\newblock Representation of divisors on hyperelliptic curves, {G}r\"{o}bner
  bases and plane curves with quasi-toric relations, 2021.
\newblock arXiv:2102.05794.

\bibitem{Toku10}
H.-o. Tokunaga.
\newblock Geometry of irreducible plane quartics and their quadratic residue
  conics.
\newblock {\em J. Singul.}, 2:170--190, 2010.

\bibitem{Toku}
H.-o. Tokunaga.
\newblock \textup{Sections of elliptic surfaces and {Z}ariski pairs for
  conic-line arrangements via dihedral covers}.
\newblock {\em \textit{J. Math. Soc. Japan}}, 66(2):613--640, 2014.

\end{thebibliography}

\noindent Ryosuke Masuya\\
Department of Mathematical Sciences, Graduate School of Science,\\
Tokyo Metropolitan University, 1-1 Minami-Ohsawa, Hachiohji 192-0397 JAPAN\\
E-mail: {\tt masuya-ryosuke@ed.tmu.ac.jp}

\end{document}